\documentclass[11pt]{amsart}

\usepackage{amsmath,amsthm,amscd} 
\usepackage{paralist}
\usepackage{mathrsfs}

\usepackage[abbrev]{amsrefs}

\newtheorem{theorem}{Theorem}[section]

\newtheorem{conj}{Conjecture}
\newtheorem{lemma}{Lemma}[section]
\newtheorem{cor}{Corollary}[section]
\newtheorem{prop}{Proposition}[section] 

\newtheorem{definition}{Definition}[section]

\theoremstyle{remark}
\newtheorem{remark}{Remark}[section]

\renewcommand{\bar}{\overline}

\newcommand{\pa}{\partial}

\newcommand{\vphi}{\phi}
  
\newfont{\fnt}{cmr10 scaled 550}
\newcommand{\eps}{\varepsilon}  

\newcommand{\bu}{{\bf u}}

   \newcommand{\C}{\mathbb C}

\newcommand{\R}{\mathbb R}

\newcommand{\Z}{\mathbb Z}

\newcommand{\frk}[1]{{\mathfrak{#1}}}

\newcommand{\bb}{{\frac{\sqrt{-1}}{2\pi}}}

\newcommand{\ka}{K\"ahler }
\newcommand{\ii}{\sqrt{-1}}

\newcommand{\gd}{\delta}

\begin{document}
\title[Gauss-Bonnet-Chern theorem on moduli space]{Gauss-Bonnet-Chern  theorem on moduli space}
\date{November 14, 2008}

\author{Zhiqin Lu}
\author{Michael R.  Douglas}

 \subjclass[2000]{Primary: 53A30;
Secondary: 32C16}
\keywords{Chern classes, Calabi-Yau moduli, Weil-Petersson
metric, Hodge metric}

\email[Zhiqin Lu, Department of Mathematics, UC Irvine, Irvine, CA 92697]{zlu@uci.edu}

\email[Michael R. Douglas, Simons Center for Geometry and Physics, Stony Brook University, Stony Brook, NY 11794]{mdouglas@scgp.stonybrook.edu}

\thanks{The first  author is  supported by the NSF grant
 DMS-12-06748, and the second  author is supported by DOE grant DE-FG02-96ER40959}
\maketitle \tableofcontents \pagestyle{myheadings}

\section{Introduction}
The moduli space of complex structures of a polarized K\"ahler manifold 
is of fundamental interest to algebraic geometers and to string theorists.
The study of its geometry is greatly enriched by considering the
Hodge structure of the underlying manifolds.  The moduli space at
infinity is particularly interesting because it is related to the
degeneration of \ka manifolds.  Even when we know very little about
the degeneration itself, by the works of Schmid~\cite{schmid},
Steenbrink~\cite{steenbrink}, Cattani-Kaplan-Schmid~\cite{cks1}, and
many others, we have a good understanding of the degeneration of the
corresponding Hodge structures.

The purposes of this paper are twofold: first, based on the work
of~\cites{schmid,steenbrink,cks1}, we prove a Gauss-Bonnet-Chern type theorem
in full generality for the Chern-Weil forms of Hodge bundles.  That
is, the Chern-Weil forms compute the corresponding Chern classes.
This settles a long standing problem. Second, we apply the result to
Calabi-Yau moduli, and proved the corresponding Gauss-Bonnet-Chern type
theorem in the setting of Weil-Petersson geometry. As an application
of our results in string theory, we prove that the number of flux
vacua of type II string compactified on a Calabi-Yau manifold is
finite, and their number is bounded by an intrinsic geometric
quantity.

Several partial results  in this direction have been
proved during the last 30 years. If the dimension of the moduli space
is $1$, then Chern-Weil forms are ``good'' in the sense of
Mumford~\cite{mumford1} and thus the Gauss-Bonnet-Chern theorem
follows. Cattani-Kaplan-Schmid~\cite{cks1} and Koll\'ar~\cite{kollar}
independently obtained the Gauss-Bonnet-Chern type theorem for the first
Chern class on a general moduli space. In the setting of
Weil-Petersson geometry~\cite{ls-1}, the rationality of the
Weil-Petersson volume was proved independently in Lu-Sun~\cite{ls-2}
and Todorov~\cite{todo-1}, and the rationality of the first Chern class of the
Weil-Petersson metric was proved in~\cite{ls-2}. If the dimension of
the moduli space is $2$, then the Gauss-Bonnet-Chern type theorem was proved
in Douglas-Lu~\cite{dl-2}.

Since in general, moduli spaces are not compact, the major difficulty
in the proof of a Gauss-Bonnet-Chern type theorem is to give a good estimate
of the Chern-Weil forms at infinity. If the growth of the connections
and curvatures {\it were} mild, or in the terminology of
Mumford~\cite{mumford1}, the forms are ``good'', then by a general
theorem in~\cite{mumford1}, the Gauss-Bonnet-Chern type theorem is
valid. However, Chern-Weil forms are not ``good'' even when the
invariant polynomial is linear (corresponding to the case of the first
Chern class). In this simplest case some tricks have been used to
prove the results (cf. Koll\'ar~\cite{kollar} and Lu-Sun~\cite{ls-2}). The proof of the general case 
is far beyond the techniques developed in the above two papers. Unlike
the compact case, the ``Splitting Principle'' for characteristic
classes is not valid for moduli space.

The first main result of this paper is Theorem~\ref{thm41}.  Since it  is quite  technical, it might be a good idea to outline the proof  here.

For the sake of simplicity,  we may assume that the moduli space $M$ is smooth, and we are dealing with a single Hodge bundle. After passing a finite cover of the manifold, we  proved the following result:
\[
\int_M f(\bb R)\in \Z,
\]
where $f$ is an invariant polynomial with integer coefficients, and  $R$ is the curvature of the  Hodge bundle.

By the classical results of Mumford and Viehweg (see~\cite{v2} for details), we know that 
the moduli space must be a quasi-projective variety. Let $\bar M$ be the compactification of $M$. We  assume that $\bar M$ is smooth and the divisor $Y=\bar M\backslash M$ is of normal crossings. By the nilpotent orbit theorem, after passing a finite cover (See Lemma~\ref{lem41} for details),   the Hodge bundle extends to $\bar M$.

Let $\Gamma^0$ be a smooth connection on the Hodge bundle on $\bar M$, then we  should prove 
that
\begin{equation}\label{p3q4}
\int_M f(\bb R)=\int_{\bar M} f(\bb R^0),
\end{equation}
where $R^0$ is the curvature of $\Gamma^0$ ({Throughout this paper, we shall use  notations like $\Gamma,  R$ to represent the connection or curvature operators as well as the corresponding  matrices under a fixed  frame}). The right hand side of the above equation is an integer by the Gauss-Bonnet-Chern theorem on compact manifold, and thus the theorem follows.

We did prove~\eqref{p3q4}. However, the proof is not straightforward. If~\eqref{p3q4} is true for one smooth connection $\Gamma^0$, it must be true for any smooth connection. The complicated asymptotic behaviors of the Hodge metric indicated by the $SL_2$-orbit theorem reveals that we have to construct a smooth connection $\Gamma^0$ very carefully in order that  it matches different growth rates of the Hodge metric in different directions.

To take a closer look at this  phenomenon, we let $\tilde f$ be the polarization of $f$. There are two ways to compare the Chern-Weil forms $f(R)$ and $f(R^0)$: the standard  way is to write $f(R)-f(R^0)$ in terms of the Bott-Chern form. 
The estimation of the Bott-Chern form at infinity of moduli space should be interesting by itself but we were not able to get the appropriate results.
This might not be very surprising because of the lack of intrinsic background metrics.

The other way is more elementary.  By the linearity of the polarization $\tilde f$, we get 
\[
\tilde f(R,\cdots,R)-\tilde f(R^0,\cdots,R^0)=\sum_j
\tilde f(\underbrace{R,\cdots,R}_j,\bar\pa(\Gamma-\Gamma^0),R^0,\cdots,R^0).
\]
Using  integration by parts, 
we were able to reduce the proof of the theorem to the proof of the uniform (with respect to $\eps$) Poincar\'e boundedness (Definition~\ref{wet}) of  the form
\begin{equation}\label{p4q3}
\bar\pa\rho_\eps\wedge\tilde  f({R,\cdots,R},\Gamma-\Gamma^0,R^0,\cdots,R^0),
\end{equation}
where $\rho_\eps$ is the cut-off function defined in Lemma~\ref{lem42} and $\Gamma$ is the connection of the Hodge metric.
The proof of  the Poincar\'e boundedness of ~\eqref{p4q3}  is the most difficult part of the paper.

The   estimate of ~\eqref{p4q3} being local, we consider the differential form on the space $\Delta^{*a}\times \Delta^b$, where $\Delta$ and $\Delta^*$ are the unit disk and the punctured unit disk, respectively. 
Let's assume that the coordinates on $\Delta^{*a}\times \Delta^b$ are $(s_1,\cdots,s_a,w_1,\cdots,w_b)$. The set $\{s_1\cdots s_a=0\}$ is called the {\it infinity} of $\Delta^{*a}\times \Delta^b$. Obviously, $\Gamma$ and $R$ are not smooth at infinity. What we know about $\Gamma$ and $R$ is from
the $SL_2$-orbit theorem~\cite{cks1} of several variables: there are finitely many cones $C_j$ with $\Delta^{*a}\times \Delta^b=\bigcup C_j$ such that, after a singular change of frame, the connection and curvature  matrices are Poincar\'e bounded.

More precisely, let
 the singular frame change be represented by  the transition matrix ${\bf e_j}$ (which blows up at infinity). 
Then
$Ad({\bf e_j}) \Gamma$ and $Ad({\bf e_j}) R$ are Poincar\'e bounded. The problem is that,
 since ${\bf e_j}$ blows up at infinity, even though  both $\Gamma^0$ and $R^0$ are smooth,
 $Ad({\bf e_j}) \Gamma^0$ or $Ad({\bf e_j}) R^0$ may  not be Poincar\'e bounded.
 
We overcame the difficulty by constructing a special smooth connection such that both  $Ad({\bf e_j}) \Gamma^0$ and  $Ad({\bf e_j}) R^0$ are Poincar\'e bounded. Lemma~\ref{47} is the technical heart of the paper addressing this key property. With the construction, by using the invariance of the function $\tilde f$, we get
\[
\tilde f(R,\cdots,\Gamma-\Gamma^0,\cdots,R^0)=
\tilde f(Ad({\bf e_j}) R,\cdots,Ad({\bf e_j})(\Gamma-\Gamma^0),\cdots,Ad({\bf e_j})R^0),
\]
which is Poincar\'e bounded and this finished the proof of the theorem.

The second major result of this paper is Theorem~\ref{thmm52},
in which we prove the convexity property of Chern-Weil forms and give an intrinsic bound of the forms in terms of the generalized Hodge metrics.

The Hodge metric on classifying  space was studied by Griffiths and Schmid as early as in ~\cite{gs} (and later by Peters ~\cite{peters}). The curvature properties (with respect to the Hermitian connection) of the Hodge metric were the key to prove the nilpotent and $SL_2$-orbit theorems. 
In~\cite{Lu3}, the first  author took a further step and proved that the Hodge metric, when restricted to a horizontal slice, must be K\"ahlerian. Moreover, the \ka metric has negative holomorphic sectional curvatures\footnote{The result is different from that of Griffiths-Schmid, because the negativity of the curvature of a {\it submanifold}, instead of the whole classifying space, was proved.} as well as other good curvature properties. By  the Schwarz lemma of Yau~\cite{Y3}, the Hodge metric must be Poincar\'e bounded.

In~\cite{fl1}, the results in~\cite{Lu3} were generalized to degenerate cases, and the generalized Hodge metrics were defined. Again, the generalized Hodge metrics are Poincar\'e bounded.

It would  be appropriate to make a remark on the Schwarz lemma at this moment.
The one dimensional Schwarz lemma (the Kobayashi's hyperbolicity) was used extensively in~\cite{schmid} and is one of the most important technical tool in that paper. Although the high dimensional Schwarz lemma is not absolutely necessary in the proof of the $SL_2$-orbit theorem of multi-variables, we did need the lemma to get  intrinsic estimate. In particular, Theorem~\ref{thmm52} is not only intrinsic, but also sharper than the estimate in~\cite{cks1}*{(5.22)}.

The differential geometry on a Calabi-Yau moduli space is called the
Weil-Petersson geometry~\cite{ls-1}. It is important in both complex
geometry and string theory.  On any smooth part of a moduli space, the
Weil-Petersson metric can be defined. But on Calabi-Yau moduli space,
the Weil-Petersson metric is directly related to the variation of
Hodge structure, thanks to a theorem of Tian~\cite{t1}. In
\S\ref{sec6}, we discovered the relation between our results on Hodge
bundles to those in Weil-Petersson geometry.  The main result of that
section, Theorem~\ref{thm63}, essentially follows from
Theorem~\ref{thm41} and Theorem~\ref{thmm52}. However, in
Weil-Petersson geometry, we can do more. We feel that Theorem~\ref{thmuj} has  no Hodge theoretic counterpart. There should be one more
layer of convexity in the Weil-Petersson geometry
(Conjecture~\ref{conj1}). By the computation in~\cite{lu12}, the
conjecture can be interpreted as the property of the underlying
Calabi-Yau manifolds. We will include results in this direction in a
separate paper.

Our results on Weil-Petersson geometry have many applications in
string theory.  The geometry of Calabi-Yau moduli and its singularity
structure is important in mirror symmetry \cite{morrison-4} and in
deriving field theoretic limits of string compactifications
\cite{kachru-vafa}.

In \cite{HorneMoore}, Horne and Moore suggested that finiteness of the volume
of a moduli space of string compactifications should be important in
early cosmology, so that realistic vacua could be produced with
finite probability.

In~\cite{AD}, Ashok and Douglas began the project of counting the
number of flux compactifications of the type IIb string on a
Calabi-Yau threefold.  If our Universe is described by string theory,
this is one factor in a count of the number of feasible Universes,
which in the type II context appears to be dominant.

In~\cite{AD}*{equation (1.5)}, a formula for the index counting these
supersymmetric vacua was given. Moreover, in Corollary 1.6 and Theorem
1.8 of ~\cite{DSZ-3}, an estimate for the error term of formula (1.5)
of Ashok-Douglas~\cite{AD} was obtained. By our results, all these
numbers are finite. Thus by the theory of Ashok and Douglas, the
number of these parallel Universes is finite.  Note also
\cite{eguchi-tachikawa} which makes similar observations in a number
of special cases.

{\bf Acknowledgement.} 
During the last four years, we have discussed the topic with many mathematicians and string theorists. We thank them all. 
The first author particularly thanks J. Koll\'ar and W. Schmid for stimulating  discussions during the preparation of this paper.

\section{Variation of  Hodge structures}\label{sec2}

Most of the material in this section is standard. It can be found in ~\cites{Gr, gh, gs}.

Let $Z$ be a compact K\"ahler manifold
of dimension $n$. A $(1,1)$ form  $\omega$ is called a polarization 
if $[\omega]$ is the first Chern class of an ample line bundle over $Z$.
The pair 
 $(Z,\omega)$, or sometimes $Z$ itself,  is
called a  polarized \ka manifold.

Let $H^k(X,\C)$ be the $k$-th cohomology group of $Z$ with $\C$ as the coefficient. 
Using $\omega$, we can define
\[
L: H^k(Z,\C)\rightarrow H^{k+2}(Z,\C), 
\quad
[\alpha]\mapsto [\alpha\wedge\omega]
\]
to be the multiplication by $[\omega]$ for $k=0,\cdots, 2n-2$. Let $P^k(Z,\C)={\rm ker}\, L^{n-k+1}$ 
on $H^k(Z,\C)$ for $k\leq n$. The group $P^k(Z,\C)$  is  called the $k$-th primitive cohomology group.

Let $H^k(Z,\Z)$ be the $k$-th cohomology group with $\Z$ as the coefficient. Let $H^{p,q}(Z)$ be the complex space of harmonic $(p,q)$ forms.
In view of  the Hodge decompositon theorem and the two Lefschetz theorems, we make the following definition:
let  $H_\Z=P^k(Z,\C)\cap H^k(Z,\mathbb Z)$ and  
$H^{p,q}=P^k(Z,\C)\cap H^{p,q}(Z)$ for $p+q=k$. Let $H_{\R}=H_\Z\otimes \R$ and $H=H_\Z\otimes\C$.
Suppose that $S$ is the quadratic form on $H_\Z$ 
induced by the cup product of the cohomology group
$H^k(X,\C)$.  
$S$ can be represented by
\begin{equation}\label{polar1}
S(\phi, \psi)=(-1)^{k(k-1)/2}\int_Z \phi\wedge \psi\wedge\omega^{n-k}
\end{equation}
for $\phi, \psi\in H$.
By the following definition, the triple $\{H_\Z, H^{p,q},S\}$ defines the so-called polarized Hodge structure.

\begin{definition}
A polarized Hodge structure
of weight $k$, which is denoted by $\{H_\Z, H^{p,q},S\}$, or
$\{H_\Z,F^p,S\}$, is given by a lattice $H_\Z$, a decomposition
\[
H=\oplus H^{p,q}
\]
with $p+q=k$ and 
\[
H^{p,q}=\overline{H^{q,p}},
\]
together with a bilinear form
\begin{equation}\label{2-1}
S: H_\Z\otimes H_\Z\rightarrow \Z,
\end{equation}
which is skew-symmetric if $k$ is odd and symmetric if $k$ is 
even.  The bilinear form satisfies the two Hodge-Riemann
relations:

\begin{enumerate}\label{hr}
\item  $S(H^{p,q},H^{p',q'})=0 \quad unless\quad p'=k-p,
q'=k-q$;
\item $(\sqrt{-1})^{p-q}\,S(\phi,\bar\phi)>0$ {\it for any
nonzero  element} $\phi\in H^{p,q}$.
\end{enumerate}
Alternatively, the decomposition~\eqref{2-1} can be described by the
filtration $\{F^p\}$ of $H$:
\[
0\subset F^k\subset F^{k-1}\subset\cdots\subset F^0= H,
\]
such that
\[
H=F^p\oplus\bar{F^{k-p+1}}, \quad H^{p,q}=F^p\cap\overline{F^q}.
\]
In this case, the Hodge-Riemann relations can be written as
\begin{enumerate}
\item[(3)]  $S(F^p,F^{k-p+1})=0$ for $p=1,\cdots, k$;
\item[(4)]  $(\sqrt{-1})^{p-q}\,S(\phi,\bar\phi)>0$  {\it for any
nonzero  element} $\phi\in H^{p,q}$.
\end{enumerate}
\end{definition}

\begin{definition}\label{def22}
The dual classifying space $\hat D$ for the 
polarized Hodge structure of weight $k$ is the set of all filtrations
\[
0\subset F^k\subset\cdots\subset F^1\subset F^0=H, \qquad
F^p\oplus\bar{F^{k-p+1}}=H,
\]
or the set of all the decompositions
\[
\sum H^{p,q}=H,\qquad H^{p,q}=\bar{H^{q,p}}\qquad (p+q=k)
\]
on which $S$ satisfies the 
Hodge-Riemann relation 1 or 3 above. The classifying space $D$ is an open set of $\hat D$ defined by the Hodge-Riemann relation 2 or 4 above.
\end{definition}

$D$ and $\hat D$ are dual in the following Lie group theoretic sense:
let
\begin{equation}\label{df23}
G_\R=\{\xi\in {\rm Hom}(H_\R,H_\R)\mid S(\xi \phi,\xi\psi)=S(\phi,\psi)\}.
\end{equation}
Then $D$ can also be written as the homogeneous space 
\begin{equation}\label{df24}
D=G_\R/V,
\end{equation}
where $V$ is the compact subgroup of $G_\R$ which leaves a fixed Hodge
decomposition $\{H^{p,q}\}$ invariant. $G_\R$ is a
semisimple real Lie group of noncompact type without compact factors.

Let $M$ be the compact dual of $G_\R$, and let $G_\C$ be the complexification of $G_\R$. Then $V\subset M\cap G_\R$ and 
\[
\hat D=M/V=G_\C/B,
\]
where $B$ is some parabolic subgroup of $G_\C$. 

\smallskip

Over the classifying space $D$ we have the holomorphic vector bundles
$\underline{F}^k,
\cdots, \underline{F}^1, \underline{H}$ whose fibers at each point are  $F^k,\cdots, F^1,H$, respectively. These bundles are called
Hodge bundles.

We   identify the holomorphic tangent
bundle
$T^{(1,0)}(D)$  as a subbundle of $Hom(\underline
H,\underline H)$:
\begin{equation}\label{k-s-2}
T^{(1,0)}(D)\subset \oplus Hom(\underline F^p,\underline H/\underline
F^p)=
\underset{r> 0}{\oplus}
Hom( \underline H^{p,q},\underline H^{p-r,q+r}),
\end{equation}
such that the following compatibility condition holds
\[
\begin{CD}
F^{p} @> >>F^{p-1}\\
@VVV        @VVV\\
H/F^p@ <<< H/F^{p-1}
\end{CD}.
\]

\begin{definition}  
A subbundle $T_h(D)$  is called the
horizontal distribution of $D$, if
\[
T_h(D)=\{\xi\in T^{(1,0)}(D)\mid  \xi F^p\subset F^{p-1}, p=
1,\cdots, k\}.
\]
\end{definition}

\begin{definition}
A horizontal slice ${\mathcal M}$ of $D$ is a complex integral submanifold
of the distribution $T_h(D)$. A holomorphic map $\phi: {\mathcal M}\rightarrow D$ is called horizontal, if its image is a horizontal slice at the regular values.
\end{definition}

 $T_h(D)$ is not always integrable. Thus in general, $\dim\mathcal M\leq\,{\rm rank}\, T_h(D)$.

Let $U$ be an open neighborhood of the universal deformation 
space (Kuranishi space)
of a polarized \ka manifold $Z$. Assume that $U$ is smooth. Then for each $Z'$ near $Z$, we
have the  isomorphism $H^*(Z',\C)=H^*(Z,\C)$ induced from the differmorphism between $Z$ and $Z'$. Under this isomorphism, 
$\{H^{p,q}(Z')\cap P^k(Z',\C)\}_{p+q=k}$ can be regarded as a
point of
$D$. The map 
\[
\vphi: U\rightarrow D,\qquad Z'\mapsto \{H^{p,q}(Z')\cap
P^k(Z',\C)\}_{p+q=k}
\]
is called the period map. 

In general, the universal deformation space allows some singularities so that the monodromy group
$\Upsilon$  is not trivial.  Thus the period map is actually a map from $U$ to $\Upsilon\backslash D$. Let $\tilde U$ be the universal covering space of $U$. Then $\vphi$ lifts to $\tilde \vphi:\tilde U\rightarrow D$. We will call both $\vphi$ and $\tilde \vphi$ period maps.

The most important property of the period map is the following~\cite{Gr}:

\begin{theorem}[Griffiths]\label{thm21}
The period map $\vphi: U\rightarrow \Upsilon\backslash D$ or $\tilde \vphi: \tilde U\rightarrow D$ is holomorphic. Furthermore, it is  an
immersion and  
 is 
 horizontal.
\end{theorem}

\section{Asymptotic behavior of the period map}
We don't know much about the degeneration of a general given family of compact K\"ahler manifolds.
But  the asymptotic behaviors of the corresponding Hodge structures were known to Schmid~\cite{schmid},
Steenbrink~\cite{steenbrink}, and Cattani-Kaplan-Schmid~\cite{cks1} by their  results of nilpotent and $SL_2$-orbit theorems of one and several variables. These theorems are not only deep, but also very long. For the full version of the theorems, we refer to the above papers. In this section, we will only define and discuss what we need for the rest of the paper. Most of the materials of this section can be found in~\cites{schmid, ck2,cks1}.

We first introduce the nilpotent orbit theorem of several variables.
 Let
$f: {\mathcal X}\rightarrow \mathscr U$ be a family of compact polarized \ka manifolds.
Since the  study of  degeneration of   Hodge
structures is local, we assume that $\mathscr U=\Delta^{*a}\times \Delta^{b}$,
where $\Delta$, $\Delta^*$ are the unit
disk and the punctured unit disk in the complex plane, respectively.
Let $(s_1,\cdots,s_a,w_1,\cdots,w_b)$ be the standard coordinate system of $\mathscr U$.
Consider the  period map
\[
\phi:\Delta^{*a}\times \Delta^{b}\rightarrow \Upsilon\backslash D,
\]
where $\Upsilon$ is the monodromy group. Let $V$ be the upper half plane. Then
$V^a\times \Delta^{b}$ is the universal covering space of $\Delta^{*a}\times \Delta^{b}$ , and we can lift $\phi$ to a map
\[
\tilde\phi: V^a\times \Delta^{b}\rightarrow D.
\]
Let $(z_1,\cdots,z_a,w_1,\cdots,w_b)$ be the coordinates of $V^a\times \Delta^b$ such that  $s_j=e^{2\pi i z_j}$ for $1\leq j\leq a$.
Corresponding to each of the first $a$ variables, we choose a
monodromy transformation $T_j\in \Upsilon$,  so that
\begin{align*}
&\tilde\phi(z_1,\cdots, z_{j}+1,\cdots, z_a,w_{1},\cdots,w_{b})\\
&\qquad=T_j\circ\tilde\phi(z_1,\cdots,z_a,w_{1},\cdots,w_{b})
\end{align*}
holds identically in all variables. $T_j$'s commute
with each other.  By a theorem of Borel,
after passing a finite cover, we may assume that the eigenvalues of $T_j$ are all $1$ (we call such $T_j$'s {\it unipotent}),
so that we can define
$N_j=\log T_j$ using the Taylor expansion of the logarithmic   function.  These $N_j$'s are called nilpotent operators.  All $N_j$'s    commute  with each other. 

Let ${\bf z}=(z_1,\cdots,z_a)$,  and
${\bf w}=(w_{1},\cdots,w_{b})$. The map
\begin{equation}\label{nil-1}
\tilde\psi({\bf z},{\bf w})=\exp(-\sum_{j=1}^az_jN_j)\circ\tilde\phi({\bf z},{\bf w})
\end{equation}
remains invariant under the translation $z_j\mapsto z_j+1, 1\leq
j\leq a$. It follows that $\tilde\psi$ drops to a map
\[
\psi:\Delta^{*a}\times\Delta^{b}\rightarrow {D}\hookrightarrow\hat D.
\]

\begin{theorem}[Nilpotent Orbit Theorem~\cite{schmid}]\label{tnot}
The map
$\psi$ extends holomorphically to $\Delta^{a+b}$. For ${\bf w}\in
\Delta^{b}$, the point
\[
F({\bf w})=\psi(0,{\bf w})\in\hat D
\]
is left fixed by $T_{j}, 1\leq j\leq a$. 
For any given number $\eta$ with $0<\eta<1$, there exist constants
$\alpha,\beta\geq 0$, such that under the restrictions
\[
{\rm Im}\, z_j\geq\alpha, 1\leq j\leq a\quad{\rm and}\quad
|w_j|\leq\eta, 1\leq j\leq b,
\]
the point $\exp(\sum_{j=1}^az_jN_j)\cdot F({\bf w})$ lies in $D$ and
satisfies the inequality~\footnote{The inequality is a refinement of ~\cite{schmid}*{(4.12)}, which was  observed by Deligne. See ~\cite{cks1}*{page 465} for details.}
\[
d(\exp(\sum_{j=1}^az_jN_j)\cdot F({\bf w}), \tilde\phi({\bf z},{\bf w}))
\leq
C\sum_{j=1}^a({\rm Im\,} z_j)^\beta  \exp(-2\pi {\rm Im \,}z_j),
\]
where $C$ is a constant and $d$ is the $G_\R$ invariant Riemannian distance function on
$D$. Finally, the mapping
\[
({\bf z},{\bf w})\mapsto \exp(\sum_{j=1}^az_jN_j)\cdot F({\bf w})
\]
is horizontal. 
\end{theorem} 

\qed

In complying with the above theorem, we make the following definition
\begin{definition}
A nilpotent orbit is a map $\theta: \C^a\rightarrow \hat D$ of the form
\begin{equation}\label{nil-2}
\theta({\bf z})=\exp (\sum_{j=1}^a z_j N_j)\cdot F,
\end{equation}
where

\begin{compactenum}[i.)]
\item $F\in\hat D$;
\item $\{N_j\}^a_{j=1}$ is a commuting set of nilpotent elements of $\mathfrak g_\R$, the Lie algebra of $G_\R$;
\item $\theta$ is horizontal, that is, $N_j(F^p)\subset F^{p-1}$;
\item There exists an $\alpha\in\R$ such that $\theta({\bf z})\in D$ for ${\rm Im\,}(z_j)>\alpha$.
\end{compactenum}
\end{definition}

\smallskip

\smallskip

Given a (real) nilpotent endomorphism $N$ of a finite dimensional complex vector space $H$, we  consider
 the monodromy weight filtration $W=W(N)$. This is defined as the unique increasing filtration $(W_j)_{j\in \Z}$ (defined over $\R$) satisfying
\begin{compactenum}[i.)]
\item $N(W_l)\subset W_{l-2}$;
\item For every $l\geq 0$, $N^l: Gr_l^W\rightarrow Gr_{-l}^W$ is an isomorphism.
\end{compactenum}
Assume that $N^{k+1}=0$. Unless otherwise stated, we will use a shifting
\[
W_l(N,k)=W_{l+k}(N).
\]
Note that $W_l(N,k)=\{0\}$ for $l<0$ and $W_l(N,k)=H$ for $l\geq 2k$.

We have the following abstract Lefschetz decomposition:
\[
Gr_l^W=\bigoplus_{j\geq 0}N^j(P_{l+2j}),
\]
where the {\it primitive} subspaces $P_l\subset Gr_l^W$ are defined by
\[
\begin{array}{ll}
P_l={\rm ker}\,\{N^{l+1}\mid Gr_l^W\rightarrow Gr_{-l-2}^W\},&{\rm if}\,\, l\geq 0,\\
P_l=\{0\},&{\rm if }\,\, l<0.
\end{array}
\]

Let $F\in\hat D$, $N\in\frk g_0$ be a nilpotent element such that $N^{k+1}=0$, and $W$ be an increasing filtration on $H$. We shall say that $(W,F,N)$, or sometimes $(W,F)$,  is a {\it polarized mixed Hodge structure}, if
\begin{compactenum}[i.)]
\item $W$ is the monodromy weight filtration of $(N,k)$;
\item $(W,F)$ is a mixed Hodge structure; that is, for $l\geq 0$, the filtration induced by $F$ on $Gr_l^W$ is a Hodge structure of weight $l$;
\item $N(F^p)\subset F^{p-1}$ for $0\leq p\leq k$;
\item For $l\geq k$, the Hodge structure induced by $F$ on the primitive subspace $P_l\subset Gr_l^W$ is polarized by the bilinear form $S_l$, where $S_l(\cdot,\cdot)=S(\cdot,N^l\cdot)$, and $S$ is a nondegenerate bilinear form on $H$ which is  skewsymmetric if $k$ is odd and symmetric if $k$ is even.
\end{compactenum}

\smallskip

A splitting of a mixed Hodge structure $(W,F)$ is a bigrading $H=\oplus J^{p,q}$ such that
\[
W_l=\underset{p+q\leq l}\bigoplus J^{p,q},\qquad F^p=\underset{r\geq p}\bigoplus J^{r,s}.
\]
We have
\[
J^{p,q}\equiv \bar{J^{q,p}} \,({\rm mod}\, W_{p+q-1})
\]
for any splitting of $(W,F)$.  A mixed Hodge structure is called to split  over $\R$, if $J^{p,q}=\bar{J^{q,p}}$.  An $(r,r)$-morphism $X$ of $(W,F)$ is compatible with $\{J^{p,q}\}$
 if $X(J^{p,q})\subset J^{p+r,q+r}$. To any mixed Hodge structure $(W,F)$ on $H=H_\R\otimes_\R\C$ we associated the nilpotent algebra
 \[
 L^{-1,-1}_R= L^{-1,-1}_R(W,F)=\{X\in \frk g\frk l(H)\mid X(J^{p,q})\subset\bigoplus_{r\leq p-1,s\leq q-1} J^{r,s}\}.
 \]
 
 The following result is from~\cite{cks1}*{Proposition 2.20}:

 \begin{prop}\label{hprop-1}
 Given a mixed Hodge structure $(W,F)$, there exists a unique $\delta\in L^{-1,-1}_R(W,F)$ such that  $(W, e^{-i\delta}\cdot F)$ is a mixed Hodge structure which splits over $\R$. Every morphism of $(W,F)$ commutes with $\delta$; thus, the morphisms of $(W,F)$ are precisely those morphisms of $(W, e^{-i\delta}\cdot F)$ which commute with this element.
 \end{prop}

Now we turn to the several variable cases. Let $N_1,\cdots,N_a$ be a commutative set of nilpotent operators coming from the unipotent monodromy operators $T_1,\cdots, T_a$. Let

\begin{equation}\label{cone}
C=\{\lambda_1N_1+\cdots+\lambda_aN_a\mid \lambda_j>0, 1\leq j\leq a\}
\end{equation}
be the cone of the monodromy operators. The basic fact about the above setting is the following~\cites{schmid,ck2,cks1}:

\begin{theorem}
Any $N\in C$ defines the same monodromy weight filtration. Moreover, let $F$ be defined in~\eqref{nil-2}. Then $(W, F)$ defines a mixed Hodge structure, polarized by each $N\in C$.
\end{theorem}

We  call that the mixed Hodge structure is polarized by $(N_1,\cdots,N_a)$.

\smallskip

To consider the boundary of the monodromy cone $C$ in~\eqref{cone}, we need the notation of {\it relative weight filtration}. We use the settings in~\cite{cks1}*{pp. 505}. Let $W^0$ be an increasing filtration of $H$ and let $N$ be a nilpotent endomorphism of $H$ which  preserves $W^0$. Delinge~\cite{deligne-1}*{1.6.13} has shown that there exists at most one weight  filtration $W=W(N,W^0)$ of $H$ such that
\begin{enumerate}
\item[i)] $N(W_l)\subset W_{l-2}$;
\item[ii)] For each $j,l\geq 0$
\[
N^l: Gr_{l+j}^WGr_j^{W^0}\rightarrow Gr_{-l+j}^WGr_j^{W^0}
\]
is an isomorphism.
\end{enumerate}

For each set of indices $I=\{i_1,\cdots,i_r\}\subset \{1,\cdots,a\}$ let $C_I$ denote the cone spanned by $N_{i_1},\cdots, N_{i_r}$. All elements of $C_I$ define the same weight filtration $W(C_I)$.

One can prove (cf.~\cite{cks1}*{Proposition  4.72}) that given $I,J\subset\{1,\cdots,a\}$, $W(C_{I\cup J})$ is the weight filtration of any $T\in C_J$ relative to $W(C_I)$.

\smallskip

\smallskip

We need the result in~\cite{cks1} on the asymptotic behavior of the Hodge length. We begin with the following setting:
by ~\eqref{nil-1}, we have 
\[
\tilde \phi({\bf z},{\bf w})=\exp(\sum z_jN_j)\cdot \psi({\bf s},{\bf w}),
\]
where $s_j=e^{2\pi i z_j}$. 
Let $(W, \psi(0,0))$ be the mixed Hodge structure, polarized by $(N_1,\cdots,N_a)$. By Proposition~\ref{hprop-1}, we let
$F=\exp(-i\delta)\cdot\psi(0,0)$ be the $\R$-split mixed Hodge structure associated  to $(W,\psi(0,0))$. Then $\xi=\exp(-i\delta)\cdot \psi: \Delta^{a+b}\rightarrow \hat D$
is holomorphic and $\xi(0)=F\in\hat D$. Since the complex orthogonal group $G_\C$ of the flat form $S$ acts transitively and holomorphically on $\hat D$, we can write $$\xi({\bf s},{\bf w})=u({\bf s},{\bf w})\cdot F$$ for a $G_\C$-valued function $u({\bf s},{\bf w})$, holomorphic on $\Delta^{a+b}$ and such that $$u(0,0)=1.$$
We choose a specific lifting of $u({\bf s},{\bf w})$ as follows. The $\R$-split mixed Hodge structure $(W, F)$ determines the bigrading $H=\oplus I^{p,q}$, where  $I^{p,q}=F^p\cap \bar{F}^q\cap W_{p+q}$.  Let the corresponding bigrading of the Lie algebra be  $\frk g=\oplus\, \frk g^{p,q}$. Because $F^p=\oplus_{r\geq p} I^{r,s}$, the Lie algebra of the stabilizer of $F$ in $G_\C$ is $\frk g_F=\oplus_{p\geq 0}\frk g^{p,q}$ and has the nilpotent subalgebra $\frk v=\oplus_{p<0}\frk g^{p,q}$ as a linear complement. A standard argument shows that the map $V\mapsto \exp V\cdot F$ is a holomorphic  diffeomorphism from $\frk v$ onto a neighborhood of $F$ in $\hat D$. We can write 
\[
\xi({\bf s},{\bf w})=\exp V({\bf s},{\bf w})\cdot F
\]
 for a unique holomorphic $V: \Delta^{a+b}\rightarrow\frk v$. The function $\gamma({\bf z},{\bf w})=\exp i\delta\cdot\exp\sum z_jN_j\cdot\exp V({\bf s},{\bf w})$ takes values in the unipotent subgroup $\exp \frk v$ and we have
\begin{equation}\label{stan}
\tilde\phi({\bf z},{\bf w})=\gamma({\bf z},{\bf w})\cdot F,\quad {\bf z}\in U^a, {\bf w}\in\mathscr C,
\end{equation}
where $\mathscr C$ is a small neighborhood of $\Delta^b$ at the origin. Let $C_\phi$ be the Weil operator. That is, $C_\phi=(\ii)^{p-q}$ on $H^{p,q}$. Let $W^{(j)}=W(N_1,\cdots,N_j)\,(1\leq j\leq a)$ be the weight filtration of $(N_1,\cdots,N_j)$. Define
\[
H\cong \oplus Gr_1^W(H),\qquad Gr_1^W(H)\overset{def}=Gr_{l_a}^{W^{(a)}}\cdots Gr_{l_1}^{W^{(1)}}(H).
\]

We recall the following theorem in~\cite{cks1}*{Theorem 5.21}:

\begin{theorem}\label{thm33}
 Let $(N_1,\cdots,N_a)$ be the nilpotent operators of  a variation of polarized Hodge structures over $\Delta^{*a}\times \Delta^b$, given with a specific ordering. If $v\in \bigcap_j W^{(j)}_{l_j}, \text{ and }Gr_1^W(v)\neq 0$, then (and only then)
\[
||v||\sim \left(\frac{\log |s_1|}{\log|s_2|}\right)^{l_1/2}\left(\frac{\log |s_2|}{\log|s_3|}\right)^{l_2/2}\cdots
(-\log |s_a|)^{l_a/2}
\]
on any region of the form
\[
\left\{ ({\bf s},{\bf w})\in \Delta^{*a}\times\Delta^b\,\left|\,\,\frac{\log |s_1|}{\log|s_2|}>\eps,\cdots,-\log|s_a|>\eps, {\bf w}\in\mathscr C\right.\right\}
\]
 for any $\eps>0$ and  $\mathscr C\subset \Delta^b$ compact, where
 \[
 ||v||^2=S(C_\phi \gamma({\bf z},{\bf w})\cdot v, \gamma({\bf z},{\bf w})\cdot \bar v).
 \]
\end{theorem}

\qed

 \section{The rationality of  Chern-Weil forms on moduli space}

 Let $\mathcal M$ be the moduli space of a polarized \ka manifold. By the period map $\phi$ defined in \S\ref{sec2}, the Hodge bundles $\underline F^k,\underline H$, and $\underline H^{p,q}$ can be pulled back to holomorphic bundles $\underline{\mathcal F}^k,\underline{\mathcal H}$ and $\underline{\mathcal H}^{p,q}$ on $\mathcal M$, respectively. These bundles are Hermitian vector  bundles with respect to the polarization $S(\cdot,\cdot)$ and the Weil operator $C_\phi$.  For the sake of convenience, we still call them (and the bundles $\underline{\mathcal F}^{p,q}$ defined below) Hodge bundles.

 Let $p<q$. Consider the holomorphic bundle
\[
\underline{\mathcal F}^{p,q}=\underline {\mathcal F}^p/\underline{\mathcal  F}^q.
\]
Note that 
if $q=p+1$, then $\underline{\mathcal F}^{p,q}=\underline{\mathcal H}^{p,k-p}$.

Let $\pi_{p,q}$ be the orthogonal projection operator on $\mathcal {\underline F}^{p,q}$. That is, let $\Omega\in\underline {\mathcal F}^p$ be a local holomorphic section. Then $\pi_{p,q}\Omega\in \underline {\mathcal F}^p$ and $S(\pi_{p,q}\Omega,\bar\Omega')=0$ for any $\Omega'\in\underline {\mathcal F}^q$. $\pi_{p,q}\Omega$ defines a holomorphic section of the bundle $\underline {\mathcal F}^p/\underline{\mathcal  F}^q$, but in general, it is not a holomorphic section of the bundle $\underline {\mathcal F}^p$. By abusing the notations, we will use  $\pi_{p,q}\Omega$  as  the section of both $\underline {\mathcal F}^p$ and $\underline {\mathcal F}^p/\underline{\mathcal  F}^q$.

Let $\Omega, \Omega'$ be  local holomorphic sections of $\underline {\mathcal F}^p$. Then they define  holomorphic sections of $\underline{\mathcal F}^{p,q}$. The Hodge  metric\footnote{See the footnote on page~\pageref{page-25}.} of $\underline{\mathcal F}^{p,q}$ is defined as
\begin{equation}\label{hodgemetric}
\langle \Omega,\bar\Omega'\rangle=S(\pi_{p,q}C_\phi\Omega,\overline{\Omega'}).
\end{equation}

Let $M$ be a quasi-projective submanifold of $\mathcal M$. For the sake of simplicity, we use $\underline{\mathcal F}^k,\underline{\mathcal H}$, $\underline{\mathcal H}^{p,q}$, and $\underline{\mathcal F}^{p,q}$ for both the 
bundles on $\mathcal M$ and their restrictions on $M$.

The main result of this section is:
 
 \begin{theorem}\label{thm41}
Let $M$ be a quasi-projective subvariety  of the moduli space ${\mathcal M}$ of a polarized \ka manifold. Let $R_{p,q}$ be the curvature tensor of $\underline {\mathcal F}^{p,q}$ with respect to the metric ~\eqref{hodgemetric}. Let $f_{p,q}$ be an invariant polynomial of ${\rm Hom}(\underline {\mathcal F}^{p,q},\underline {\mathcal F}^{p,q})$ with rational coefficients. Then for any sequence 
\[
(p_1,q_1),\cdots,(p_r,q_r),
 \]
of nonnegative integers, we have
\begin{equation}\label{concl}
\int_{M} f_{p_1,q_1}(\bb R_{p_1,q_1})\wedge\cdots\wedge
 f_{p_r,q_r}(\bb R_{p_r,q_r})
\in\mathbb Q.
\end{equation}
Moreover, the Chern-Weil form $\prod f_{p_j,q_j}(\bb R_{p_j,q_j})$ computes the corresponding Chern class on any compactification of $M$.
 \end{theorem}
 
 \begin{remark}
 The curvatures of the Hodge bundles blow up at infinity of the moduli space. Thus the precise meaning of ~\eqref{concl}, {\it at this time},  is 
  \[
\underset{\eps\rightarrow 0}{\lim} \int_{M} \rho_\eps f_{p_1,q_1}(\bb R_{p_1,q_1})\wedge\cdots\wedge
 f_{p_r,q_r}(\bb R_{p_r,q_r})
\in\mathbb Q,
\]
where $\rho=\rho_\eps$ is the cut-off function with compact support in $M$ such that $0\leq\rho\leq 1$, and $\rho\equiv 1$ outside the $\eps$-neighborhood of the infinity. However, in the next section, we shall prove that the integration converges absolutely so that the result is independent of the choice of the cut-off functions and the expression of ~\eqref{concl} makes sense.
 \end{remark}
 
  Let $M_{reg}$
 be the smooth part of $M$. 
 By Hironaka's theorem, there is a projective manifold $\bar M$, called the {\it compactification} of $M$, such that
 $\bar M \backslash M_{reg}$ is a divisor of $\bar M$ of normal crossings\footnote{It may be more accurate to say that, up to a finite cover, $\bar{M}$ can be chosen to be a manifold.}. 
 In order to study the asymptotic behavior of the curvatures of Hodge bundles, we first need to extend the bundles to the compactification of $M$.
 
The following lemma is due to Kawamata~\cite{kawa3}.  For the sake of completeness, we sketch the proof here. Note that in~\cite{ls-2}, we proved the similar result under the setting  of Weil-Petersson geometry.
  
 \begin{lemma}  \label{lem41} 
 By replacing $M$ with $M_{reg}$, we assume that $M$ is smooth. Moreover, we assume that $\bar M$ is smooth and $Y=\bar M\backslash M$ is a divisor of normal crossings. Then there is a finite branched cover $\tilde M$ of 
 $\bar M$ with branched locus $Y$ such that the monodromy operators along $Y$ in $\tilde {M}$ are all unipotent.
 \end{lemma}
 
 {\bf Proof.} 
Let $T$ be a
monodromy operator along certain irreducible component of the divisor 
$Y$ which is not unipotent. Let
$T=\gamma_{s}\gamma_{u}$ be the decomposition of $T$ into its
semi-simple part and its unipotent part. By a theorem of Borel, 
there is an integer $d$ such that $\gamma_s^d=1$. Let $L$ be an
ample line bundle of $\bar M$. 
Let $Y=\sum D_j$ be the decomposition of the divisor
$Y$ into irreducible components. 
 We assume that the monodromy
operator $T$ is generated by $U\backslash D_1$,
where $U$ is a  neighborhood of $D_1$.
Assume that $t$ is large enough
such that the bundle $L^t(-D_1)$ is very ample. By taking the
$d$-th root of the defining section of $L^t(-D_1)$ we get a variety
$M_1$ such that outside a  divisor, it is a finite
covering space of $M$. $M_1$ may have some singularities. However,
we can always remove those 
divisors containing singularities to get a smooth manifold. 

Let $M'$ be an Zariski open set of $M_1$ that is a covering space of $M$.
The Hodge
bundles
can be pulled back to the  manifold $M'$. At any neighborhood
 of $\bar{M}\backslash  M$, the
transform of
$(\bar M, M)$  to $(M_1,{M'})$ is the $d$-branched cover defined by
$z_1\mapsto
\sqrt[d]{z_1}$, where $z_1=0$ locally defines  the divisor
$D_1$.  Evidently, the monodromy operator $T$ is transformed to
$T^d$, which becomes  a unipotent operator. 

We observe that if $T'$ is a unipotent  monodromy operator,
then
 under the transform 
$(\bar M, M)$ to $(M_1,{M'})$, $T'$ is
still unipotent.   
Since there are only finitely many irreducible components
of $Y$, there are only finitely many monodromy operators which are not unipotent.\footnote{During the process, it is possible that some divisors are added. However, along these divisors, the monodromy operators are the identity operator.}
Thus after finitely many  transforms, we can
get a  compact complex   manifold $\tilde M$ on which all
monodromy operators are unipotent.

\qed

For the rest of the section, we let $Y$ be the divisor of $\tilde{M}$ of normal crossings such that $\tilde{M}\backslash Y$ is a finite covering of $M$. 
We also assume that $\tilde {M}$ is covered by   finite open coordinate neighborhoods $\{U_\alpha\}_{\alpha=1,\cdots,t}$, and $\{\psi_\alpha\}$ is the partition of unity subordinating to the cover.

Using Lemma~\ref{lem41} and the nilpotent orbit theorem (Theorem~\ref{tnot}), the Hodge bundles $\underline{\mathcal F}^p$, $\underline{\mathcal H}^{p,q}$, $\underline {\mathcal H}$, and $\underline{\mathcal F}^{p,q}$ extend to vector bundles over $\tilde{M}$. We use the same notations to denote these (extended) bundles. 

The following definition is  standard to experts. However, we find the notation quite convenient to use:

\begin{definition} \label{wet}
Let $\Delta, \Delta^*$ be the unit disk and the punctured unit disk of $\mathbb C$, respectively. Let $U=\Delta^{*a}\times\Delta^b$ and let the standard coordinate system of $U$ be $(s_1,\cdots,s_a,w_1,\cdots,w_b)$. A differential form on $U$ is called  Poincar\'e  bounded, if its components are  bounded under the coframe
\[
\frac{ds_i}{s_i\log|s_i|},\quad \frac{d\bar s_i}{\bar s_i\log|s_i|}, dw_j, d\bar w_j
\]
for $i=1,\cdots,a, j=1,\cdots,b$.
\end{definition}

The following Lemma is obvious

\begin{lemma}\label{easy}
If a form is bounded, then it is Poincar\'e bounded. 
If $\sigma_1, \sigma_2$ are Poincar\'e bounded, so is $\eta_1\wedge\eta_2$. 
\end{lemma}

\qed

Moreover,  the notation of Poincar\'e boundedness is independent of the choice of coordinates:

\begin{lemma}\label{easy-two}
Let $M$ be a quasi-projective manifold and let $\bar M$ be its compactification such that $\bar M\backslash M$ is a divisor of normal crossings. Let $U,U'$ be two neighborhoods of the divisor such that $U\cap U'\neq\emptyset$. Then a smooth form  is Poincar\'e bounded on any compact subset of $U\cap U'\cap M$   with respect to the coordinate system of $U$ if and only if it is Poincar\'e bounded  with respect to  that of $U'$.
\end{lemma}

{\bf Proof.} We assume that $U\approx \Delta^{*a}\times\Delta^b$ and 
$U'\approx \Delta^{*a'}\times\Delta^{b'}$. Let the coordinates of the two neighborhoods be $(s_1,\cdots,s_a,w_1,\cdots,w_b)$ and $(s_1',\cdots,s_{a'}',w_1',\cdots,w_{b'}')$, respectively. We further assume that the divisor is the zero locus of either $\{s_1\cdots s_a=0\}$ or $\{s_1'\cdots s_{a'}'=0\}$. Assuming $a\leq a'$, then on $U\cap U'$, we can rearrange the order of $s_j'$ such that
\[
s_j=\xi_j s_j'
\]
for $j=1,\cdots,a$, where $\xi_j$ are smooth nonzero functions. Since
\[
d\log s_j=d\log s_j'+d\log \xi_j,
\]
we concluded that 
$\frac{ds_j}{s_j\log|s_j|}$, $\quad \frac{d\bar s_j}{\bar s_j\log|s_j|}$ are bounded under the coframe $\frac{ds_j'}{s_j'\log|s_j'|}$, $\quad \frac{d\bar s_j'}{\bar s_j'\log|s_j'|}$ for $j=1,\cdots,a$, and vice vesa. On the other hand, $s_j'$ for $j>a$ are 
nonzero.  The lemma is proved.

\qed

As the first step of the proof of Theorem~\ref{thm41}, we need to choose the cut-off function on $\tilde{M}$ carefully. The following construction of the cut-off function depends on the particular choice of the cover and  partition of unity. However, the main feature  is that the complex Hessian of the function is of order $O(\frac{1}{\eps^2(\log\frac 1\eps)^2})$, a little bit better than $O(\frac{1}{\eps^2})$. 

\begin{lemma}\label{lem42}
For any $\eps>0$ small enough, there is a smooth real function $\rho=\rho_\eps$ on $\tilde{M}$ such that
\begin{enumerate}
\item $0\leq\rho\leq 1$;
\item $\pa\rho, \bar\pa\rho$, and $\pa\bar\pa\rho$ are Poincar\'e bounded;
\item The Euclidean measure of ${\rm supp}\, (\pa\rho)$ goes to zero as $\eps\rightarrow 0$;
\item In a neighborhood of $Y$, $\rho\equiv 0$;  and $\rho(x_0)= 1$ if the distance of $x_0\in M$ to $Y$ is greater than $2\eps$.
\end{enumerate}
\end{lemma}

{\bf Proof.} 
Rearranging the order of $\{U_\alpha\}$, we may assume that
$\{U_1,\cdots, U_s\}$ is an open cover of the divisor $Y$ and  $$(U_{s+1}\cup\cdots \cup U_t)\cap Y=\emptyset.$$
 We further assume that $U_\alpha\backslash Y=(\Delta^*)^{a_\alpha}\times \Delta^{b_\alpha}$ and the coordinates are 
$(s_1^\alpha,\cdots, s^\alpha_{a_\alpha}, w_1^\alpha,\cdots,w^\alpha_{b_\alpha})$. Assume that locally $Y$ is the zero locus of 
\[
s_1^\alpha\cdots s_{a_\alpha}^\alpha=0
\]
on each $U_\alpha$.
Let $\eta: \mathbb R\rightarrow \mathbb R$, $0\leq \eta\leq 1$ be a smooth decreasing function defined as 
\[
\eta(x)=\left\{
\begin{array}{ll} 
0 & x\geq 1\\
1 & x\leq 0
\end{array}
\right..
\]

Let 
\begin{equation}\label{phi-1}
\eta_\eps(z)=\eta\left(\frac{(\log\frac 1 {|z|})^{-1}-\eps}{\eps}\right),
\end{equation}
and let
\[
\eta_\eps^\alpha(s_1^\alpha,\cdots,s_{a_\alpha}^\alpha)=\prod_{j=1}^{a_\alpha} (1-\eta_\eps (s_j^\alpha)).
\]
Then the function $\rho=\rho_\eps$ is defined as
\[
\rho_\eps=\sum_{\alpha=1}^s\psi_\alpha\eta_\eps^\alpha+\sum_{\alpha=s+1}^t\psi_\alpha,
\]
where $\{\psi_\alpha\}$ is the fixed partition of unity defined before.

For the above $\rho_\eps$, (1) is trivial. In order to prove (2), by Lemma~\ref{easy} and~\ref{easy-two}, we only need to prove that $\bar\pa\eta_\eps$ and $\pa\bar\pa\eta_\eps$ are Poincar\'e bounded. By a straightforward computation, we have
\begin{align*}
&\bar\pa\eta_\eps=\frac 1{2\eps}\eta'\frac{d\bar z}{\bar z(\log\frac{1}{|z|})^2};\\
&\pa\bar\pa\eta_\eps=\frac 1{4\eps^2}\eta''\frac{dz\wedge d\bar z}{|z|^2(\log\frac{1}{|z|})^4}
+\frac 1{2\eps}\eta'\frac{dz\wedge d\bar z}{|z|^2(\log\frac{1}{|z|})^3}.
\end{align*}
The above expressions are non-zero only if
\[
\eps\leq\left(\log\frac{1}{|z|}\right)^{-1}\leq 2\eps.
\]
Thus there is a constant $C$ such that 
\[
|\bar\pa\eta_\eps|\leq C\left|\frac{ d\bar z}{|z|(\log\frac{1}{|z|})}\right|,\quad
|\pa\bar\pa\eta_\eps|\leq C\left|\frac{dz\wedge d\bar z}{|z|^2(\log\frac{1}{|z|})^2}\right|,
\]
and thus both $\bar\pa\eta_\eps$ and $\pa\bar\pa\eta_\eps$ are Poincar\'e bounded.

Since (3) is implied by (4), we only prove the latter. Let $x_0\in M$. For any $x_0$ which is close enough to $Y$, $\psi_\alpha=0$ for $\alpha\geq s+1$. On the other hand, since $x_0$  is close to $Y$, for any $1\leq\alpha\leq s$, there is an $s^\alpha_{j(\alpha)}$ which is sufficiently small. Thus $\eta_\eps(s_{j(\alpha)}^\alpha)=1$ and consequently $\eta_\eps^\alpha=0$. This proves $\rho_\eps(x_0)=0$. If the distance of $x_0$ to $Y$ is at least 
$2\eps$, then there is a constant $C>0$ such that $|s_j^\alpha|\geq C\eps$ for any $1\leq j\leq a_\alpha$ and $1\leq\alpha\leq s$. Since $\eps\log\eps^{-1}\to 0$ for $\eps$ small, we have $\eta_\eps^\alpha=1$ for $1\leq\alpha\leq s$. Thus we conclude that $\rho_\eps(x_0)=\sum\psi_\alpha(x_0)=1$ and this completes the proof.

\qed

By pulling back the Hodge bundles to $\tilde{M}\backslash Y$,  the curvature operator $R_{p_j,q_j}$ makes sense as the ${\rm Hom}\, (\underline{\mathcal F}^{p_j,q_j},\underline{\mathcal F}^{p_j,q_j})$-valued
 $(1,1)$-form on $\tilde{M}\backslash Y$. 
Since  $\tilde{M}\backslash Y$ is a finite cover of $M$, there is a positive integer $\mu$ such that \begin{align*}&
 \int_{\tilde{M}} f_{p_1,q_1}(\bb R_{p_1,q_1})\wedge\cdots\wedge
 f_{p_r,q_r}(\bb R_{p_r,q_r})\\&=\mu\int_{M} f_{p_1,q_1}(\bb R_{p_1,q_1})\wedge\cdots\wedge
 f_{p_r,q_r}(\bb R_{p_r,q_r}).
 \end{align*}
 
 Let $c_1^j,\cdots, c_{d_j}^j$ be the elementary invariant polynomials on ${\rm Hom}\,(\mathbb C^{d_j},\mathbb C^{d_j})$, where
 $d_j$ is the rank of the vector bundle $\underline{\mathcal F}^{p_j,q_j}$. 
  Then there are polynomials $g_j$ such that
 \[
 f_{p_j,q_j}=g_j(c_1^j,\cdots,c_{d_j}^j)
 \]
 for $j=1,\cdots,r$.

The following theorem implies the main result of this section, Theorem~\ref{thm41}: 

  \begin{theorem}\label{422}
   Using the same assumptions and notations as in Theorem~\ref{thm41} and Lemma~\ref{lem42}, we have
  \begin{equation}\label{42}\underset{\eps\rightarrow 0}{\lim}
   \int_{\tilde{M}} \rho_\eps f_{p_1,q_1}(\bb R_{p_1,q_1})\wedge\cdots\wedge
 f_{p_r,q_r}(\bb R_{p_r,q_r})\in\mathbb Z,
 \end{equation}
 if the coefficients of the polynomials $g_j$ ($j=1,\cdots, r$) are integers.
 \end{theorem}

Let $n_j$ be the degree of $f_{p_j,q_j}$.
 Let $\tilde f_{p_j,q_j}$ be the polarization of $f_{p_j,q_j}$. That is,
\[
\tilde f_{p_j,q_j}: (\C^{d_j\times d_j})^{n_j}\rightarrow \C,
\]
such that
\begin{enumerate}
\item $\tilde f_{p_j,q_j}(A_1,\cdots,A_{n_j})$ is linear with each $A_l$ ($1\leq l\leq n_j$);
\item $\tilde f_{p_j,q_j}$ is symmetric. That is,  $$\tilde f_{p_j,q_j}(\cdots,A_k,\cdots,A_l,\cdots)=\tilde f_{p_j,q_j}(\cdots,A_l,\cdots,A_k,\cdots);$$
\item $\tilde f_{p_j,q_j}(A,\cdots,A)=f_{p_j,q_j}(A)$.
\end{enumerate}

We let 
\[
\mathcal E=\bigoplus_{j=1}^r\underline{\mathcal F}^{p_j,q_j}.
\]
Then
\[
R=\begin{pmatrix}
R_{p_1,q_1}\\&\ddots\\&&R_{p_r,q_r}
\end{pmatrix}
\]
is the curvature tensor of $\mathcal E$. 
For the sake of simplicity, we use $R_j=R_{p_j,q_j}$, $f_j=f_{p_j,q_j}$, and  $\tilde f_j=\tilde f_{p_j,q_j}$ for $1\leq j\leq r$. Let
\[
f(R)=f_1(R_1)\wedge\cdots\wedge f_r(R_r).
\]

Let $\Gamma^0$ be a smooth $(1,0)$-type connection of $\mathcal E$ over $\tilde{M}$ of  the form
\[
\Gamma^0=\begin{pmatrix}
\Gamma^0_1\\&\ddots\\&&\Gamma^0_r
\end{pmatrix},
\]
where $\Gamma^0_j$ are smooth $(1,0)$-type connections on $\underline{\mathcal F}^{p_j,q_j}$. Let
$R^0=\bar\pa \Gamma^0$ be the curvature tensor. Since $\tilde M$ is compact, we have
\begin{equation}\label{gb}
\int_{\tilde{M}} f(\bb R^0)\in\mathbb Z
\end{equation}
by the Gauss-Bonnet-Chern theorem on smooth manifold.

We proved the following

\begin{lemma}\label{lem44}
If there is a smooth connection $\Gamma^0$ on $\tilde{M}$ such that for any $i$ and $j$, $f_j(R_j)$ and 
\[
\eta_{i,j}=
\bar\pa\rho_\eps\wedge \tilde f_j(\underbrace{R_j,\cdots,R_j}_i, \Gamma_j-\Gamma_j^0, R_j^0,\cdots,R_j^0)
 \]
 are Poincar\'e bounded, then Theorem~\ref{422} is true. Here $\Gamma_j$ is  the connection operator of the Hodge bundle $\underline{\mathcal F}^{p_j,q_j}$.
 \end{lemma}

{\bf Proof.}  We use the following obvious equality
\begin{align*}
& f_1(R_1)\wedge\cdots\wedge f_r(R_r)-f_1(R^0_1)\wedge\cdots\wedge f_r(R^0_r)\\
&=\sum_{j=1}^r\left\{\left(\prod_{l<j} f_l(R_l)\right)\wedge(f_j(R_j)-f_j(R_j^0))\wedge 
\left(\prod_{l>j}f_l(R_l^0)\right)\right\}.
\end{align*}
Noting that $\Gamma_j-\Gamma_j^0$ is globally defined for any $j$, we have
\begin{align}\label{leb}
\begin{split}
&\underset{\eps\rightarrow 0}{\lim}\int_{\tilde{M}}
\rho_\eps( f(R)-f(R^0))\\
&=-\underset{\eps\rightarrow 0}{\lim}\int_{\tilde{M}}
\sum_{i,j}\eta_{i,j}\wedge
\left(\prod_{l<j} f_l(R_l)\right)\wedge \left(\prod_{l>j}f_l(R_l^0)\right).
\end{split}
\end{align}
By the assumption and by Lemma~\ref{easy}, the integrand of the right hand side of the above equation is  Poincar\'e bounded. Let $\xi ds_1\wedge\cdots \wedge d\bar s_a\wedge dw_1\wedge \cdots \wedge d\bar w_b$ be the $(a+b,a+b)$-component of the integrand on a general neighborhood $U=U_\alpha$, where $\xi$ is a smooth function on $U$. Then there is a constant $C$ such that
\[
|\xi|\leq C\prod_{j=1}^a\frac{1}{|s_j|^2(\log\frac{1}{|s_j|})^2}.
\]
It is elementary to see that the above function is Euclidean integrable. Since the Euclidean measure of ${\rm supp}\, (\pa\rho_\eps)$ goes to zero, by the Lebesgue theorem, 
the right hand side of ~\eqref{leb} is zero.

The lemma and hence Theorem~\ref{thm41}  follow from~\eqref{gb}, the ordinary Gauss-Bonnet-Chern Theorem.

\qed

Before giving the explicit construction of the connection $\Gamma^0$, we define the local frame on each $U=U_\alpha$ ($1\leq \alpha\leq s$). We use the notations in  ~\eqref{stan}. Let $F_\infty^p$ be the limiting Hodge filtration. Then for any basis $\{\tilde v_{p,j}\}$ 
of $F_\infty^p$,
\[
\exp (\sqrt{-1}\delta)\exp V({\bf s},{\bf w}) \tilde v_{p,j}
\]
gives a local frame of the bundle $\underline{\mathcal F}^p$. In fact, this is the local frame we use to define the extension of the Hodge bundles. We call such a local frame {\it defined by the nilpotent orbit theorem}. Likewise, if $\{v_{p,q,j}\}$ is a basis of the vector space
\[
\bigoplus_{j=1}^r F^{p_j}_\infty/F^{q_j}_\infty,
\]
then 
\[
\exp (\sqrt{-1}\delta)\exp V({\bf s},{\bf w}) v_{p,q,j}
\]
gives a local frame of the bundle $\mathcal E$, which we also call it defined by the nilpotent orbit theorem.

Now we construct the connection $\Gamma^0$ explicitly: as before $U_\alpha\cap Y\neq\emptyset$ if and only if $1\leq \alpha\leq s$. On each $\underline{\mathcal F}^{p_j,q_j}$, if $1\leq\alpha\leq s$, let $\Omega_{\alpha,p_j,q_j,a}$, where $a=1,\cdots,d_j$, be the local holomorphic frame of the bundle $\underline{\mathcal F}^{p_j,q_j}$ defined by the nilpotent orbit theorem; if $\alpha\geq s+1$, we let $\Omega_{\alpha,p_j,q_j,a}$ be an arbitrary holomorphic local frame of $\underline{\mathcal F}^{p_j,q_j}$. Let
\[
\Omega_\alpha=(\Omega_{\alpha,p_1,q_1,1},\cdots,\Omega_{\alpha,p_1,q_1,d_1},\cdots,\Omega_{\alpha,p_r,q_r,1},\cdots,\Omega_{\alpha,p_r,q_r,d_r}).
\]
Then the transition matrices of the vector bundles $A_{\alpha\beta}$ are holomorphic matrix-valued functions:
\[
\Omega_\alpha=\Omega_\beta\, {A}_{\alpha\beta}^t
\]
on $U_\alpha\cap U_\beta\neq \emptyset$.  Note that $\{A_{\alpha\beta}\}$ are (block) diagonalized matrix-valued functions. The diagonalization is compatible with respect to the direct sum structure of $\mathcal E=\bigoplus_{j=1}^r\underline{\mathcal F}^{p_j,q_j}$.

We define the connection matrix
\begin{equation}\label{p-3-q}
\Gamma_\alpha^0=\sum_\gamma\psi_\gamma\pa A_{\alpha\gamma}A_{\alpha\gamma}^{-1}.
\end{equation}
on $U_\alpha$, where $\{\psi_\alpha\}$ is the partition of unity subordinating to the cover $\{U_\alpha\}$.
Then $\Gamma_\alpha^0$ and the curvature matrix $R_\alpha^0$ are all (block) diagonalized.

As a general fact, we have the following result: 

\begin{lemma} The collection of matrix-valued $(1,0)$ forms $\{\Gamma_\alpha^0\}$ defines a smooth $(1,0)$ connection on the vector bundle $\mathcal E\rightarrow\tilde{ M}$.
\end{lemma}

{\bf Proof.} The compatibility conditions of the transition matrices are
\[
A_{\alpha\gamma}=A_{\alpha\beta}A_{\beta\gamma}
\]
on $U_\alpha\cap U_\beta\cap U_\gamma\neq\emptyset$.
Thus
\[
\pa A_{\alpha\gamma}=\pa A_{\alpha\beta}A_{\beta\gamma}+A_{\alpha\beta}\pa A_{\beta\gamma}.
\]
It follows that 
\[
\pa A_{\alpha\gamma}A^{-1}_{\alpha\gamma}=\pa A_{\alpha\beta}A^{-1}_{\alpha\beta}
+A_{\alpha\beta}\pa A_{\beta\gamma}A^{-1}_{\beta\gamma}A_{\alpha\beta}^{-1}.
\]
Using~\eqref{p-3-q}, we have
\[
\Gamma_\alpha^0=\pa A_{\alpha\beta}A^{-1}_{\alpha\beta}+A_{\alpha\beta}\Gamma_\beta^0A_{\alpha\beta}^{-1}
\]
on $U_\alpha\cap U_\beta\neq\emptyset$. Thus $\{\Gamma_\alpha^0\}$ defines a smooth connection of $\mathcal E$.

\qed

The problem to verify the assumptions in Lemma~\ref{lem44} is purely local. So we will concentrate on a typical neighborhood $U=U_\alpha$ and suppress the subscript $\alpha$ for the sake of simplicity. We assume that $U\approx \Delta^{*a}\times \Delta^b$. Also, we will use $\Gamma,\Gamma^0,R,R^0$ to represent the   connection and  curvature matrices, respectively, under the local frame defined by the nilpotent orbit theorem.

In order to study the properties of the smooth connection $\Gamma^0$, we use the following notations in~\cite{cks1}*{\S 5}\footnote{
Since in this paper we use a large set of notations, for simplicity, we
will try to keep them the same  as in~\cite{cks1}. Thus we have to sacrifice the uniqueness of the notations. For example, $\alpha$ is used in $U=U_\alpha$ and also as the subscripts of the set $I$, etc. It should be clear from the context.}.

As before,
we assume that the local coordinate system of $U$ is $$(s_1,\cdots,s_a,w_1,\cdots,w_b),$$ and $U\cap Y$ is the zero locus of $s_1\cdots s_a=0$.

 We define a cone $\mathcal C$ in $\Delta^{*a}$ by
 \[
 \mathcal C=\{|s_1|\leq|s_2|\leq\cdots\leq|s_a|\mid (s_1,\cdots,s_a)\in \Delta^{*a}\}.
 \]
 Such a cone gives an ordering of $\{1,\cdots,a\}$.
 Let 
 \[
 I=\{i_\alpha\}, \,1\leq i_1<\cdots<i_r= a.
 \]
 Define $|I|=r$, and 
 
   \begin{align}
 \begin{split}
 &t_\alpha=y_{i_\alpha}/y_{i_{\alpha+1}}\quad {\rm if }\, \alpha<r, \,\,t_r=y_a,\\
&u_\alpha^j=y_j/y_{i_\alpha}\qquad{\rm for }\, \,i_{\alpha-1}<j<i_\alpha, \\
&X_\alpha=X_\alpha(\bu)=N_{i_\alpha}+\underset{i_{\alpha-1}<j<i_\alpha}{\sum} u^j_\alpha N_j,
\end{split}
\end{align}
where $N_1,\cdots,N_a$ are the nilpotent operators.
Let $Y_\alpha({\bf u})$ be the semisimple elements corresponding to $X_\alpha$ by the Jacobson-Morozov theorem. Define
\begin{equation}\label{pe-24}
e({\bf y})=e({\bf t},{\bf u})=\exp(\sum_\alpha\frac 12 \log y_{i_\alpha}Y_\alpha({\bf u}))=\exp(\sum_\alpha\frac 12 \log t_{\alpha}{\bf Y}_\alpha({\bf u})),
\end{equation}
where ${\bf Y}_\alpha=Y_1+\cdots+Y_\alpha$.

Moreover, let 
 \begin{itemize}
 \item $\mathscr A$: analytic functions of ${\bf u}\in \mathbb R_+^{a-r}$;
 \item $\mathscr L$: Laurent polynomials in $\{t_\alpha^{1/2}\}$ (or $\{y_{i_\alpha}^{1/2}\}$ with coefficients in $\mathscr A)$;
 \item $\mathcal O$: pull back to $V^a$ of the ring of holomorphic germs at $0$ in $\Delta^a$, via $V^a\rightarrow\Delta^{*a}\rightarrow \Delta^a$, where $V$ is the upper half plane and $V\to\Delta^*$ is the covering map $z\mapsto e^{2\pi\sqrt{-1}z}$;
 \item ${\mathscr L}^{b}$: polynomials in $\{t_\alpha^{-1/2}\}$ with coefficients in $\mathscr A$;
 \item $(\mathcal O\otimes \mathscr L)^b$: subring of $(\mathcal O\otimes \mathscr L)$ generated by $\mathcal O$, ${\mathscr L}^b$, and all polynomials of the form $s_jt^{m_1/2}_{\alpha_1}\cdots t^{m_p/2}_{\alpha_p}$ for $j\leq i_{\alpha_l}$, $m_l\in\mathbb Z$, $l=1,\cdots,p$;
 \item $(\mathcal O\otimes \bar{\mathcal O}\otimes\mathscr L)^b$: ring generated by
 $(\mathcal O\otimes\mathscr L)^b$,  $(\bar{\mathcal O}\otimes\mathscr L)^b$;
 \item ${\mathscr R}^b_{K,L}$: ring of rational expressions $f/g$, $f,g\in (\mathcal O\otimes \bar{\mathcal O}\otimes\mathscr L)^b$, with $g$ bounded away from zero on $(V^a)^I_{K,L}$ (defined below), where $0<K<L$ are two positive numbers.
 \end{itemize}
Via ${\bf y}\rightarrow ({\bf t}, {\bf u})$, $\mathbb R^a_+$ is identified with $\mathbb R_+^r\times \mathbb R_+^{a-r}$. For any $K\ll L$, 
let
 \begin{align}
 \begin{split}
 & (\R_+^a)_{K,L}^I=\left\{ {\bf y}\in \R^a_+\mid t_\alpha>L, \, 1\leq u_\alpha^j\leq K \right\},\\
 & (V^a)^I_{K,L}=\left\{ {\bf z}\in V^a\mid {\bf z}={\bf x}+\sqrt{-1} {\bf y},\, {\bf y}\in (\R_+^a)^I_{K,L}\right\},\\
& (\Delta^{*a})^I_{K,L}=\left\{ {\bf s}\in \Delta^{*a}\,\left|\, \frac{\log s_j}{2\pi\sqrt{-1}}\in (V^a)^I_{K,L}\right\}\right..
 \end{split}
 \end{align}
 
 \smallskip

 \begin{remark}
 The above definition is slightly different from the one in ~\cite{cks1}*{page 509}. We use the notation $(\Delta^{*a})^I_{K,L}$ instead of $(\Delta^{*a})^I_{K}$ so that  it is a little easier to show the $C^\infty$ convergence when $L\gg K$ in  Theorem~\ref{thm33}. Of course, we do have the $C^\infty$ convergence on $(\Delta^{*a})^I_{K,K}=(\Delta^{*a})^I_{K}$. But the proof is well hidden in~\cite{cks1}*{\S 5} and more explanations are needed.   \end{remark}
 
With the above settings, we have the following combinatorial  lemma:
 
 \begin{lemma}
Let
 \[
 1=K_{a+1}\leq K_a<K_{a-1}<\cdots<K_0=+\infty
 \]
 be a sequence such that
 \[
 K_j>K_{j+1}^a
 \]
 for $j=1,\cdots,a$. Let
 \[
 A_j=\bigcup_{|I|=j}(\Delta^{*a})^I_{K^a_{j+1},K_j}
 \]
 for $j=1,\cdots,a$, and let
 \[
 A_0=\{{\bf s}\mid y_a\leq K_1+1  \}.
 \]
  Then we have
 \[
\bigcup_{j=0}^a A_j\supset\mathcal C.
 \]
 \end{lemma}

{\bf Proof.} We let
\[
\xi_1=\frac{y_1}{y_2},\,\,\dots,\,\,\xi_{a-1}=\frac{y_{a-1}}{y_a},\,\,\xi_a=y_a.
\]
Consider $(a+1)$ intervals
\[
(K_{a+1},K_a],\,\,\cdots,(K_2,K_1],\,\,(K_1,K_0(=+\infty)).
\]
By the pigeonhole principle,  if ${\bf s}\not\in A_0$, then there is an $1< l\leq a+1$ such that
\[
\xi_j\not\in (K_l,K_{l-1}]
\]
for any $1\leq j\leq a$.  Define $I=\{i_1<i_2<\cdots<i_r=a\}$ such that
\[
\begin{array}{ll}
\xi_j\leq K_l& \text{ for } j\not\in I;\\
\xi_j> K_{l-1}&\text{ for } j\in I.
\end{array}
\]
Note that we must have $i_r=a$.
Then we have
\[
u_\alpha^j=\frac{y_j}{y_{i_\alpha}}=\frac{y_j}{y_{j+1}}\cdots\frac{y_{i_\alpha-1}}{y_{i_\alpha}}
\leq K^n_l,
\]
and
\[
t_\alpha=y_{i\alpha}/y_{i_{\alpha+1}}>K_{l-1}.
\]
Thus we have ${\bf s}\in A_{l-1}$ and the lemma is proved.

\qed

In what follows, we will use matrix notations extensively. In particular, for a fixed frame, the Hodge metrics (locally) are represented by matrices, and the change of frames are represented using the matrix notations as well.

Let $\Omega=\Omega_\alpha$ be the local frame of $\mathcal E$ defined by the nilpotent orbit theorem. Let $C$ be a fixed cone of $U$ of the form $(\Delta^{*a})^I_{K,L}\times\mathscr C$, where $\mathscr C$ is a compact subset of $\Delta^b$.
In~\cite{cks1}*{pp. 514}, for such a cone, there is a basis of $\{v_{C,j}\}$ of $\bigoplus_{j=1}^r F^p_\infty/F^q_\infty$ on the typical fiber $H$ flagged according to the ``limiting split'' Hodge filtration $F$, which we call it {\it the  basis} of the cone.  We
let $\Omega_C$ be the  frame of the cone defined by the above basis via the nilpotent orbit theorem, and let ${\bf e}$ be the matrix of $e=e({\bf y})$ of ~\eqref{pe-24} under the frame $\Omega_C$. Then the following is true (cf. (5.19) or ``Proof of (5.22)'' of ~\cite{cks1}):

\begin{theorem}\label{out}
 Let ${\bf h}$ be the metric matrix of the Hodge bundle $\mathcal E$ under the  basis $\Omega_C$, and let
\[
{\bf h}={\bf e}^t{\bf k}\bar{\bf e}.
\]
Then the matrix ${\bf k}$ and its inverse matrix are bounded on the cone $C$.
\end{theorem}

\qed

By the definition of the  basis on the cone,  there is a {\it constant} matrix $A_C$ (cf.~\cite{cks1}*{(5.20)})  such that
\[
\Omega=\Omega_CA_C^t.
\]
Let $\Gamma_C$, $R_C$ be the connection and the curvature operators of the Hodge metric under the local frame $\Omega_C$, respectively. Then in Proposition (5.22) of ~\cite{cks1},   the following was proved

\begin{theorem} \label{thm445}The coefficients of the forms $Ad({({\bf e}^{-1}})^t)\Gamma_C$ and $Ad({({\bf e}^{-1}})^t)R_C$ are Poincar\'e bounded.
\end{theorem}

\qed

The key technical lemma of this section is the following:
\begin{lemma}\label{47} Let $U'=U_\gamma$ be an open set such that $U'\cap C\neq\emptyset$. Let $A=A_{\alpha\gamma}$. Then on $U'\cap C\neq\emptyset$,
\[
Ad(({\bf e}^{-1})^t) Ad(A_C^{-1}) (\pa A A^{-1})
\]
is Poincar\'e bounded.
\end{lemma}

{\bf Proof.}  Let $\Omega'=\Omega_\gamma$ be the local frame of $U'$ defined by the nilpotent orbit theorem. Let $C'$ be a fixed cone of $U'$ and let $\Omega_{C'}'$ be the  frame of the cone. We assume that $C\cap C'\neq\emptyset$. 
We just need to prove the assertion of the lemma on $C\cap C'$ because as $C'$ is running over all the cones, the whole $U'\cap C$ will be covered.

 Let ${\bf e}'$ be the matrix under the frame $\Omega'_{C'}$. Let $A_{C'}$ be the constant matrix defined as 
\[
\Omega'=\Omega'_{C'}A_{C'}^t.
\]
Then we have
\[
\Omega_C=\Omega_{C'}'A_{C'}^tA^t{(A_C^{-1})^t}.
\]
We let $B=A_C^{-1}AA_{C'}$ and let ${\bf h}'$ be the metric matrix of $\Omega_{C'}'$. Then $\Omega_C=\Omega_{C'}'B^t$. Thus
\begin{equation}\label{fol-1}
{\bf h}=\Omega_C^t\bar{\Omega_C}=B{\bf h'}{\bar B}^t.
\end{equation}
It follows that
\begin{equation}\label{fol}
\pa{\bf h}{\bf h}^{-1}=\pa BB^{-1}+Ad(B)(\pa{\bf h}'({\bf h}')^{-1}).
\end{equation}
Since $A_C, A_{C'}$ are constant matrices, we have
\[
\pa BB^{-1}=Ad(A_C^{-1})(\pa AA^{-1}).
\]
By Theorem~\ref{thm445}, from ~\eqref{fol}, we see that in order to prove the lemma, we only need to prove that
\[
Ad(({\bf e}^{-1})^t) Ad(B)\pa {\bf h}'({\bf h}')^{-1}=D (Ad((({\bf e'})^{-1})^t))(\pa {\bf h'}({\bf h'})^{-1} )D^{-1}
\]
is Poincar\'e bounded, where $D=({\bf e}^{-1})^t B({\bf e'})^t$ . Using Theorem~\ref{thm445} again, we only need to prove that $D$ and $D^{-1}$  are  bounded. To prove this, we observe that
if
\[
{\bf h}={\bf e}^t{\bf k}\bar{\bf e},
\]
and if
\[
{\bf h}'=({\bf e'})^t{\bf k}'\bar{\bf e'},
\]
then by~\eqref{fol-1}
\[
{\bf k}=D{\bf k}'\bar D^t.
\]
By Theorem~\ref{out}, ${\bf k}$, ${\bf k}'$ and their inverse matrices  are bounded. 
Since both ${\bf k}$ and ${\bf k}'$ are positive definite,
$D$ and $D^{-1}$ must be bounded.

\qed

\begin{cor}\label{cor41}
 Let $\Gamma_C^0=
Ad(A_C^{-1})\Gamma^0$  and $R_C^0=Ad(A_C^{-1})R^0$ be the connection   (note that $A_C$ is a constant matrix) and the curvature matrices  under the frame $\Omega_C$. Then
$Ad(({\bf e}^{-1})^t) (\Gamma_C^0)$ and $Ad(({\bf e}^{-1})^t) (R_C^0)$ are Poincar\'e bounded.
\end{cor}

{\bf Proof.} We have
\[
\Gamma_C^0=Ad(A_C^{-1})(\sum_\gamma\psi_\gamma\pa A_{\alpha\gamma}A^{-1}_{\alpha\gamma})=\sum_\gamma\psi_\gamma Ad(A_C^{-1})(\pa A_{\alpha\gamma}A^{-1}_{\alpha\gamma}).
\]
By the above lemma, for each $\gamma$, $Ad(({\bf e}^{-1})^t)Ad(A_C^{-1})(\pa A_{\alpha\gamma}A^{-1}_{\alpha\gamma})$ is Poincar\'e bounded. Since $\{U_\alpha\}$ is a locally finite cover, the conclusion on $\Gamma_C^0$ follows. The result on $R_C^0$ follows from a similar formula:
\[
Ad(({\bf e}^{-1})^t)R_C^0=\sum_\gamma\bar\pa\psi_\gamma
Ad(({\bf e}^{-1})^t) Ad(A_C^{-1})(\pa A_{\alpha\gamma}A^{-1}_{\alpha\gamma}).
\]

\qed

{\bf Proof of Theorem~\ref{422}.} We only need to verify the assumptions in Lemma~\ref{lem44}. We mention for one more times that all the matrix-valued functions or forms we have constructed so far are (block) diagonalized with respect to the direct sum structure of the bundle $\mathcal E$. This fact allows us to suppress the index $j$
in Lemma~\ref{lem44} when we fix  the $j$.

On any cone $C$, by the invariance of $f=f_j$ and $\tilde f=\tilde f_j$, we have
\begin{equation}\label{plk-1}
f(R)=f(Ad(({\bf e}^{-1})^t)(R_C)),
\end{equation}
and
\begin{align}\label{pq}
\begin{split}&
\bar\pa\rho_\eps\wedge \tilde f(R,\cdots,R,\Gamma-\Gamma^0,R^0,\cdots,R^0)\\&
=
\bar\pa\rho_\eps \wedge\tilde  f(R_C,\cdots,R_C,\Gamma_C-\Gamma^0_C,R^0_C,\cdots,R^0_C)
\\&
= \bar\pa\rho_\eps\wedge \tilde f(Ad(({\bf e}^{-1})^t)(R_C),\cdots,Ad(({\bf e}^{-1})^t)(R_C),\\&\quad Ad(({\bf e}^{-1})^t)(\Gamma_C-\Gamma^0_C),Ad(({\bf e}^{-1})^t)(R^0_C),\cdots,Ad(({\bf e}^{-1})^t)(R^0_C)).
\end{split}
\end{align}
By Theorem~\ref{thm445}, $Ad(({\bf e}^{-1})^t)(R_C)$ and $Ad(({\bf e}^{-1})^t)(\Gamma_C)$ are Poincar\'e bounded. On the other side, by Corollary~\ref{cor41}, $Ad(({\bf e}^{-1})^t)(R_C^0)$ and $Ad(({\bf e}^{-1})^t)(\Gamma_C^0)$ are also Poincar\'e bounded.  Thus the left hand sides of~\eqref{plk-1} and ~\eqref{pq} are Poincar\'e bounded. By Lemma~\ref{lem44}, this implies Theorem~\ref{422} (hence Theorem~\ref{thm41}).

\qed

For the moduli space ${\mathcal  M}$ itself, we have

\begin{cor} Using the same notations as in Theorem~\ref{thm41}, we have
\begin{equation}\label{concl-also}
\int_{\mathcal M} f_{p_1,q_1}(\bb R_{p_1,q_1})\wedge\cdots\wedge
 f_{p_r,q_r}(\bb R_{p_r,q_r})
\in\mathbb Q.
\end{equation}
\end{cor}

{\bf Proof.}  By the theorem of Viehweg~\cite{v2}, $\mathcal M$ is a quasi-projective variety. 
The corollary follows from Theorem~\ref{thm41}.
\qed

\section{The generalized Hodge metrics}

In this section, we shall show that the form in~\eqref{concl} is absolutely integrable. The result can be proved using Proposition (5.22) in~\cite{cks1}, or by the argument in the last section. However, we provide a proof here which is elementary (avoid using the $SL_2$-orbit theorem) and sharper. More importantly, we give the intrinsic upper bound of the integrals, which can be regarded as Chern number inequalities on  moduli space.

In the first part of this section, we recall some notations and results in~\cite{fl1}.

Let $M$ be  a complex manifold of dimensional $m$. Suppose that $M$ is  the  parameter space of a family of polarized compact \ka manifolds $\pi:\mathcal X\rightarrow M$. 
 By the 
functorial property, the Hodge bundles $\underline{\mathcal  H}, \underline{\mathcal F}^p, \underline{\mathcal H}^{p,q}$, and $\underline{\mathcal F}^{p,q}$  on $M$ can be defined as the pull-back of the Hodge bundles from  the classifying spaces. The bundles can also be identified to the relative cohomology groups as follows
\[
\underline{\mathcal H}^{p,q}=PR^q\pi_{*}\Omega^p_{\mathcal X/M}, \,\,
\underline {\mathcal F}^p=\underline {\mathcal H}^{p+q,0}\oplus\cdots\oplus
\underline {\mathcal H}^{p,q}
\]
for $p,q\geq 0$, where $\Omega^p_{\mathcal X/M}$ is the sheaf of relative holomorphic $(p,0)$ forms on $\mathcal X$. In particular, $\underline{\mathcal H}=PR^k\pi_*(\mathbb C)$. 
Let $Z_x=\pi^{-1}(x)$ for $x\in M$. Assume that $\dim Z_x=n$.
The Kodaira-Spencer map $T_x({{M}}) \rightarrow
H^1(Z_x,T^{(1,0)}Z_x)$ gives the  bundle map
\[
\frac{\pa}{\pa t_i}: \underline {\mathcal H}^{p,q}\rightarrow PR^k\pi_{*}(\C)/\underline {\mathcal H}^{p,q} 
\]
for $0\leq k\leq n$
by differentiation, where $PR^k\pi_{*}(\C)$ is the primitive part of $R^k\pi_{*}(\C)$,
and $(t_1,\cdots,t_m)$ is a local holomorphic coordinate system at $x$. The map induces the  natural bundle
 map  (compare to~\eqref{k-s-2}):
\begin{equation}\label{12-1}
T^{(1,0)}({ {M}}) \rightarrow\underset{p+q=k}{\oplus} {\rm
Hom}\,(\underline {\mathcal H}^{p,q},
PR^k\pi_{*}(\C)/\underline {\mathcal H}^{p,q}).
\end{equation}

 We make the 
 following definition of the generalized Hodge metrics:\footnote{The Hodge metrics are generally referred to the natural Hermitian metrics on Hodge bundles. After~\cites{Lu3,Lu5}, the metrics in Definition~\ref{def32} are also called the Hodge metrics or the generalized Hodge metrics, because they have the similar curvature properties as the original ones on Hodge bundles.}\label{page-25}
 
 \begin{definition}\label{def32} Assume that $0\leq k\leq n$.
 Let $h_{PH^k}$ be the  pull back of the
natural Hermitian metric on the bundle $\underset{p+q=k}{\oplus}{\rm
Hom} (\underline{\mathcal H }^{p,q}\rightarrow
PR^k\pi_{*}(\C)/\underline{\mathcal H}^{p,q})$ to
$T^{(1,0)}(M)$. 
We use $\omega_{PH^k}$ to denote the
corresponding K\"ahler forms.\footnote{That is, if
$h_{PH^k}=(h_{PH^k})_{i\bar j} dt_i\otimes d\bar t_j$, then
$\omega_{PH^k}=\frac{\sqrt{-1}}{2\pi}(h_{PH^k})_{i\bar j} dt_i\wedge
d\bar t_j$.}   In complying  to the Lefschetz decomposition theorem, we define
\begin{equation}\label{12-2}
\omega_{H^k}=\omega_{PH^k}+\omega_{PH^{k-2}}+\cdots.
\end{equation}
We call both $\omega_{H^k}$ and $\omega_{PH^k}$  the
generalized Hodge metrics. 
\end{definition}

\begin{remark}\label{rk1}
The above construction is  a generalization of the
Hodge metric defined by the first  author~\cites{Lu3,Lu5}. In fact,
 it is proved in~\cite{ls-1} that
$$\omega_{PH^n}=\omega_H,$$ the latter being the Hodge metric defined in~\cite{Lu3}.  Alternatively, the generalized Hodge metrics can also be the defined as the restriction of the invariant Hermitian metrics on the corresponding classifying spaces.
\end{remark}

Because of the possible degeneration of the action~(\ref{12-1}), a
generalized Hodge metric is only semi-positive definite; hence, it
is  a pseudo-metric. Nevertheless, it enjoys the similar ``curvature"
properties of the Hodge metric. To elaborate this, we recall the following result in~\cite{fl1}*{Proposition 2.8}:

\begin{prop}
\label{prop1}
Let $c_1(E)$ be the Ricci form  of a vector bundle $E$. Then we have
\begin{equation}\label{five}
\omega _{PH^{k}}=\sum_{0\leq p\leq
k}pc_{1}(\underline{\mathcal H}^{p,k-p}),
\end{equation}
and
\begin{equation}\label{six}
\omega _{H^{k}}=\sum_{l=0}^{[\frac k2]}\sum_{p\leq k-2[\frac k2]+2l}pc_{1}(\underline{\mathcal H}^{p,k-2[\frac k2]+2l-p})
\end{equation}
for $k\leq n$.
\end{prop}

{\bf Proof.} For the sake of completeness, we include the proof here. Fixing a $k\leq n$, we let
\[
\underline{\mathcal  F}_k^p=\underline{\mathcal H}^{k,0}
\oplus \cdots\oplus \underline{\mathcal H}^{p,k-p}
\]
for $p=0,\cdots,k$. Thus  for $q=k-p$,
\[
\underline {\mathcal H}^{p,q}=\underline  {\mathcal F}_k^p/
\underline {\mathcal F}_k^{p+1}.
\]
In terms of the curvatures, we have
\begin{equation}\label{q0}
 c_1(\underline{\mathcal H}^{p,q})=c_1(
\underline {\mathcal F}_k^p) -c_1(\underline {\mathcal F}_k^{p+1}).
\end{equation}

By the Abel summation formula, we have
\begin{equation}\label{q1}
\sum_{0\leq p\leq k}pc_{1}(\underline{\mathcal H}^{p,k-p})=c_1( \underline {\mathcal F}_k^k)+\cdots+c_1( \underline {\mathcal F}_k^1)
+c_1(\underline{\mathcal F}_k^0).
\end{equation}

Each $\underline {\mathcal F}_k^p$ is a sub-bundle of the flat bundle $
\underline {\mathcal F}_k^0 =PR^k\pi_*\C$. Let $t_1,\cdots,t_m$ be the local holomorphic
coordinate of ${{M}} $ and let the bundle map
\[
\frac{\pa}{\pa t_j}: \underline  {\mathcal F}_k^p\rightarrow
 \underline {\mathcal F}_k^0/\underline {\mathcal F}_k^p, \,1\leq j\leq m
\]
be represented by
\[
\frac{\pa\Omega_\alpha}{\pa t_j}=b_{j\alpha\mu}T_\mu,
\]
where $\Omega_\alpha$ and $T_\mu$ are the basis of $
\underline {\mathcal F}_k^p$ and $ \underline {\mathcal F}_k^0/\underline{\mathcal F}_k^p$, respectively. Then the
first Chern class can be written as \footnote{This is essentially due to~\cite{gs}. See also~\cite{Gr}*{page 34}.}
\begin{equation}\label{12-3}
c_1(\underline {\mathcal F}_k^p)=\frac{\sqrt{-1}}{2\pi}
\sum_{i,j,\alpha,\mu}b_{i\alpha\mu}\bar b_{j\alpha\mu}dt_i\wedge d\bar
t_j
\end{equation}
for $0\leq p\leq k$. By the definition of the generalized Hodge metrics,  (\ref{five}) follows from ~\eqref{q1} and ~\eqref{12-3}. (\ref{six}) follows from (\ref{five}) and
the Lefschetz decomposition theorem. The proof is completed.

\qed

\begin{cor}\label{coro51}
Using the above notations, we have
\[
d\omega_{PH^k}=d\omega_{H^k}=0.
\]
In particular, if the generalized Hodge metric is positive definite, then it defines a \ka metric.
\end{cor}

\qed

Let $g_p$ be the Hodge metric on $\underline{\mathcal H}^{p,q}$.
In what follows,
we shall compute the curvatures of $g_p$ and compare them with the generalized Hodge metrics.
For the sake of simplicity, we assume that
$\underline {\mathcal  H}^{k+1,-1}=\underline {\mathcal H}^{-1,k+1}=0$ and $\underline {\mathcal F}_k^{k+1}=0$,
$\underline {\mathcal F}_k^{-1}=\underline {\mathcal F}^0$.

Fix $k\leq n$, $p\leq k$ and $q=k-p$. Let
$\{\Omega_{p,i}\}$, $i=1,\cdots, h^{p,q}={\rm rank}\,\underline {\mathcal H}^{p,q}$ be a local holomorphic
frame of $\underline {\mathcal H}^{p,q}$.

\begin{definition}\label{523}
Let $(t_1,\cdots,t_m)$ be a holomorphic local coordinate system at
a point of $M$. We define
$\nabla_\alpha\Omega_{p,i}\in H^{p-1,q+1}$
to be the projection of $\pa_\alpha\Omega_{p,i}
=\frac{\pa}{\pa t_\alpha}\Omega_{p,i}$
to $H^{p-1,q+1}$ with respect to the bilinear form
$S(\,\,,\,\,)$.
\end{definition}

For simplicity, we shall use $(\,\,,\,\,)$ in stead of the
bilinear form $S$.
With the above notation,
\begin{equation}\label{first}
(g_p)_{i\bar j}=\langle\Omega_{p,i}, \overline{\Omega_{p,j}}\rangle
=(\ii)^{p-q}(\Omega_{p,i}, \overline{\Omega_{p,j}})
\end{equation}
 is the Hermitian metric matrix of $\underline {\mathcal H}^{p,q}$
for $p=0,\cdots,k$. Using~\eqref{12-3}, we can write the generalized Hodge metric in local coordinates as follows (cf.~\cite{fl1}):

\begin{prop} \label{prop32}
For fixed $k$,
the generalized Hodge metric matrix under the local
coordinate system $(t_1,\cdots,t_m)$  can be written as 
\begin{equation}\label{second}
h_{\alpha\bar \beta}=\sum_{p=0}^k(\ii)^{p-q-2}
g_p^{i\bar
j} (\nabla_\alpha\Omega_{p,i},\overline{\nabla_\beta\Omega_{p,j}}),
\end{equation}
where $(g_p^{i\bar j})$ is the inverse of $(g_p)_{i\bar j}$.
\end{prop}

\qed

Let $R_p$ be the curvature operator of $g_p$. Then by~\cite{Gr}*{page 34}, we have
\begin{equation}\label{ade}
(R_p)_{i\bar j\gamma\bar\delta}
=(\ii)^{p-q}(\nabla_\gamma\Omega_{p,i},
\overline{ \nabla_\delta\Omega_{p,j}})
-(\ii)^{p-q} (\bar\pa_\gamma\Omega_{p,i},\overline
{\bar\pa_\delta\Omega_{p,j}}).
\end{equation}
By the Cauchy inequality, the first term of the right hand side of the above is no more than $\sqrt{h_{\gamma\bar\gamma}h_{\delta\bar\delta}}$. It can be 
proved (cf.~\cite{fl1}*{Lemma A.5})  that the operator $\nabla$ is dual to the operator  $\bar\pa$. Thus the norm of the two operators must be the same and consequently, 
the second term of the above is also no more than $\sqrt{h_{\gamma\bar\gamma}h_{\delta\bar\delta}}$. Thus we have
\begin{equation}\label{plk}
|(R_p)_{i\bar j\gamma\bar\delta}|\leq 2 \sqrt{h_{\gamma\bar\gamma}h_{\delta\bar\delta}}
\end{equation}
for any $i,j$.

\smallskip 

Now we turn to the proof of the absolute  integrability and Chern number inequalities.

\begin{theorem}\label{thmm52} Let $c_\alpha(g_p)$ be the $\alpha$-th Chern-Weil  form of the Hodge bundle $\underline {\mathcal H}^{p,q}$ with respect to the metric $g_p$. Then we have
\begin{equation}\label{sdf-1}
|c_{\alpha_1}(g_{p_1})\wedge\cdots\wedge c_{\alpha_r}(g_{p_r})|\leq
2^m\,\omega^m,
\end{equation}
where 
\[
\omega=\sum\omega_{H^k},
\]
and $\sum\alpha_j=m$.
Similarly, we have
\begin{equation}\label{sdf-2}
|c_\alpha(g_p)\wedge\omega_{PH^k}^{m-\alpha}|\leq 2^\alpha\omega_{PH^k}^m
\end{equation}
for  $k=p+q$. 
\end{theorem}

{\bf Proof.} Let $\omega_0$ be any \ka metric of $M$. Then by Corollary~\ref{coro51}, for any $\eps>0$, $\omega+\eps\omega_0$ is a \ka metric. Suppose $(t_1,\cdots,t_m)$ is a holomorphic normal coordinate system at $x_0\in M$. 
Let $(h_{\alpha\bar\beta})$ be the metric matrix of $\omega$ under this coordinate system.
Then $h_{\gamma\bar\gamma}\leq 1$ for any $\gamma$.
We first  consider a general Chern-Weil form $c_\alpha (g_p)$.
For fixed $\gamma,\delta$, by~\eqref{plk}, we have
\begin{equation}\label{plk-28}
|(R_p)_{i\bar j\gamma\bar\delta}|\leq 2\sqrt{h_{\gamma\bar\gamma}h_{\delta\bar\delta}}\leq 2.
\end{equation}
Let
\[
(R)_{i\bar j}=(R_p)_{i\bar j\gamma\bar\delta} dt_\gamma\wedge d\bar t_\delta.
\]
Then
by definition,
\[
c_\alpha(g_p)=\left(\bb\right)^\alpha\frac{(-1)^\alpha}{\alpha!}\sum_{\tau\in S_\alpha}{\rm sgn}(\tau) R_{i_1\bar{i_{\tau(1)}}}
\wedge\cdots\wedge  R_{i_\alpha\bar{i_{\tau(\alpha)}}},
\]
where $S_\alpha$ is the symmetric group on the set $\{1,2,\cdots,\alpha\}$. 
Let
\[
dt_I=dt_{i_1}\wedge\cdots\wedge dt_{i_\alpha}.
\]
Using~\eqref{plk-28},  if we write
\[
c_\alpha(g_p)=\left(\frac{1}{2\pi}\right)^\alpha\frac{1}{\alpha!}\sum a_{I\bar J} dt_I\wedge d\bar t_J,
\]
then we have
\begin{equation}\label{awdf}
|a_{I\bar J}|\leq 2^\alpha.
\end{equation}
Since the number of partition of the set $\{1,\cdots,m\}$ into subsets of $\alpha_1,\cdots,\alpha_r$ elements, respectively, is
\[
\frac{m!}{(\alpha_1!)\cdots(\alpha_r!)},
\]
using ~\eqref{awdf}, we have
\begin{equation}
\left|\frac{c_{\alpha_1}(g_{p_1})\wedge\cdots\wedge c_{\alpha_r}(g_{p_r})}{(\omega+\eps\omega_0)^m}\right|
\leq 2^m,
\end{equation}
and ~\eqref{sdf-1} follows  by taking $\eps\to 0$. Similarly, we get~\eqref{sdf-2}. The theorem is proved.

\qed

\begin{cor}[Chern number inequalities] Using the same notations as above, we have
\[
\int_M c_{\alpha_1}(g_{p_1})\wedge\cdots\wedge c_{\alpha_r}(g_{p_r})\leq
2^m\,\int_M\omega^m,
\]
and
\[
\int_M
c_\alpha(g_p)\wedge\omega_{PH^k}^{m-\alpha}\leq 2^\alpha\int_M\omega_{PH^k}^m.
\]
\end{cor}

\qed

Now we  state the following result, which can be viewed as a degenerate version of Yau's Schwarz lemma~\cite{Y3}. The result here is a slight generalization of~\cite{fl1}*{Theorem A.1} because $\tau$ doesn't have to be the Poincar\'e metric.

\begin{theorem} \label{ptheorem} Let $\tau=\bb\tau_{\alpha\bar\beta} dt_\alpha\wedge d\bar t_\beta$ be a \ka metric on $M$ such that
\begin{enumerate}
\item $\tau$ is a complete metric;
\item The Ricci curvature of $\tau$ has a lower bound.
\end{enumerate}
Then there is a constant $C$, depending only on the dimension of $M$ and the lower bound of the Ricci curvature of $\tau$, such that
\[
\omega_{H^k}\leq C\tau.
\]
\end{theorem}

{\bf Proof.}
Let $\xi$ be the smooth function defined by $\xi=\tau^{\alpha\bar\beta} h_{\alpha\bar\beta}$, where $\tau^{\alpha\bar\beta}$ is the inverse matrix of $\{\tau_{\alpha\bar\beta}\}$. Then by the Bochner type formula in ~\cite{fl1}*{Appendix A}, there is a constant $C>0$, depending only on the dimension of $M$ and the lower bound of the Ricci curvature of $\tau$, such that
\[
\Delta\xi\geq\frac{1}{C}\xi^2-C\xi.
\]
Using the generalized maximum principle~\cite{cy2} (see also~\cite{ty1}), $\xi\leq C^2$ is a bounded function. This completes the proof.

\qed

A typical choice of the metric $\tau$ is the so-called Poincar\'e metric whose  \ka form is denoted as  $\omega_P$.
 At any point $x_0\in \bar M\backslash M$, there is a neighborhood $U$ of $p$ such that $U\cap M$ can be identified as $\Delta^{*a}\times\Delta^b$. The metric $\omega_P$ on $\Delta^{*a}\times\Delta^b$ is asymptotic to the Poincar\'e metric
\[
\omega_P\sim\bb\left(\sum_{i=1}^a\frac{ds_i\wedge d\bar s_i}{|s_i|^2(\log\frac{1}{|s_i|})^2}
+\sum_{i=a+1}^{a+b}dw_i\wedge d\bar w_i\right).
\]
See ~\cite{ls-2}*{\S 5} for the detailed constructions. 

In our terminology, Theorem A.1 of~\cite{fl1} can be rephrased  as

\begin{cor}\label{cor53} The generalized Hodge metrics are Poincar\'e bounded.
In particular, the Hodge volume is finite, hence ~\eqref{concl} is absolutely integrable.
\end{cor}

{\bf Proof.} By a straightforward computation, the volume of the Poincar\'e metric  is finite. This proves the corollary.

\qed

\smallskip

The above corollary implies the result in ~\cite{cks1}*{(5.23)}
\begin{cor}  The Chern-Weil forms extend to a current on the compactification $\bar M$ of $M$.
\end{cor}

\qed

\section{Chern classes on Calabi-Yau moduli}\label{sec6}
In this section, we assume that $Z$ is a polarized Calabi-Yau manifold of dimension $n$ and $\mathcal M$ is the moduli space of $Z$ (the Calabi-Yau  moduli). 
For a Calabi-Yau manifold, the Hodge structure of weight $n$ is the most important one. 
Let $D$ be the classifying space  defined in Definition~\ref{def22} corresponding to the weight $n$, and let 
$\phi: \mathcal M\rightarrow D$ be  the period map. By~\cite{bg-1}, the map is an immersion on the smooth part of $\mathcal M$.

\begin{definition} Let $Z$ be a polarized Calabi-Yau manifold with the Ricci flat \ka metric $\mu$ whose \ka form defines the polarization. Let
$X,Y\in H^1(Z, T^{(1,0)}Z)$.  Define the $L^2$ inner product by
\[
(X,Y)=\frac{1}{n!}\int_Z\langle X,Y\rangle \mu^n.
\]
For a Calabi-Yau manifold, the Kodaira-Spencer map: $$T_Z\mathcal M\rightarrow H^1(Z, T^{(1,0)}Z)$$ is an isomorphism. Thus the above inner product defines a metric on the smooth part of $\mathcal M$. The metric happens to be K\"ahlerian, and is called the
Weil-Petersson metric of  $\mathcal M$.
\end{definition}

Let $\underline {\mathcal F}^n$ be the first Hodge bundle on $\mathcal M$.  It is a line bundle because $\dim H^{n,0}(Z)=1$ for Calabi-Yau manifolds. By Griffiths~\cite{Gr}*{page 34} and~\cite{bg-1}, we know that
\[
c_1(\underline {\mathcal F}^n)>0.
\]

The crucial result we are going to use is the following~\cite{t1}:

\begin{theorem}[Tian] 
\label{thm61} On the smooth part of $\mathcal M$, the \ka form of the Weil-Petersson metric $\omega_{WP}$
is  $c_1(\underline {\mathcal F}^n)$. More precisely, let $\Omega$ be a holomorphic local section of $\underline {\mathcal F}^n$, then
\[
\omega_{WP}=-\bb\pa\bar\pa\log(\Omega,\bar\Omega).
\]
\end{theorem}

\qed

The Weil-Petersson metric is defined on the smooth part of the Calabi-Yau moduli. Though in general, a moduli space 
may have singularities, by the theorem of Tian~\cite{t1} (see also Todorov~\cite{to}), we know that the Kuranish space of a Calabi-Yau manifold is smooth. As a result, a Calabi-Yau moduli space must be a complex orbiford, and the Weil-Petersson metric is a \ka orbifold metric.

By Theorem~\ref{thm61}, we know that the Weil-Petersson metric depends only on the variation of Hodge structure.  Moreover, we shall see that
  the curvature of the Weil-Petersson metric
is also directly related to the variation of Hodge structure, which allows us to make use of Theorem~\ref{thm41} and Theorem~\ref{thmm52} in the previous sections.

We cite the formulas of Strominger~\cite{s} (for the Calabi-Yau moduli of a Calabi-Yau threefold) and Wang~\cite{wang1} (for the general case) on the curvature tensor of the Weil-Petersson metric.

\begin{theorem}\label{thm62-1}
 On the Calabi-Yau moduli of a polarized Calabi-Yau $n$-fold,
let $\Omega$ be a local holomorphic section 
of $\underline {\mathcal F}^n$. Then  
\[
R(\omega_{WP})_{\alpha\bar\beta\gamma\bar\delta}=g_{\alpha\bar\beta} g_{\gamma\bar\delta}+g_{\alpha\bar\delta}g_{\gamma\bar\beta}
-\frac{(\nabla_\alpha\nabla_\gamma\Omega,\overline{\nabla_\beta\nabla_\delta\Omega})}
{(\Omega,\bar\Omega)}
\]
for any $1\leq\alpha,\beta,\gamma,\delta\leq m$, where $m=\dim \mathcal M$; $g_{\alpha\bar\beta}$ is the metric matrix of the Weil-Petersson metric; and $\nabla$ is the connection defined in Definition~\ref{523}.
Moreover, In the case of moduli space of Calabi-Yau threefolds, the curvature can be  represented  in terms of the Yukawa coupling
$\{F_{\alpha\beta\gamma}\}$:
\[
R(\omega_{WP})_{\alpha\bar\beta\gamma\bar\delta}=g_{\alpha\bar\beta} g_{\gamma\bar\delta}+g_{\alpha\bar\delta}g_{\gamma\bar\beta}-\sum_{\xi,\eta}\frac{1}{|(\Omega,\bar\Omega)|^2}g^{\xi\bar\eta}F_{\alpha\gamma\xi}\bar{F_{\beta\delta\eta}}.
\]
\end{theorem}

Based on the above formulas, in~\cite{ls-1}*{Theorem 5.3}, it was proved that 
the Weil-Petersson curvature is  $L^1$ (with respect to the Weil-Petersson metric).  In what follows, we shall 
greatly generalize the result. By relating the curvature tensor of the Weil-Petersson metric to  that of the Hodge bundles, we prove the corresponding Gauss-Bonnet-Chern theorem for Calabi-Yau moduli.

For the sake of simplicity, we write $R=R(\omega_{WP})$. Consider the bundle $\underline{\mathcal F}^{n-1}/\underline{\mathcal F}^n=\underline{\mathcal H}^{n-1,1}$. By~\eqref{ade}, the curvature tensor $R'$ of $\underline{\mathcal H}^{n-1,1}$ can be written as
\begin{equation}\label{1+2}
(R')_{i\bar j\gamma\bar\delta}
=(\ii)^{n-2}\left((\nabla_\gamma\Omega_{i},
\overline{ \nabla_\delta\Omega_{j}})
-(\bar\pa_\delta\Omega_{i},\overline
{\bar\pa_\gamma\Omega_{j}})\right),
\end{equation}
where $\{\Omega_1,\cdots,\Omega_n\}$ is a local frame of $\underline{\mathcal H}^{n-1,1}$.

Since $\mathcal M$ is the Calabi-Yau moduli, we have the Kodaira-Spencer isomorphism $T\mathcal M=\underline{\mathcal H}^{n-1,1}$. Let $\Omega$ be a local holomorphic section of the line bundle $\underline{\mathcal F}^n$ and let $(t_1,\cdots,t_m)$ be a holomorphic local coordinate system of $\mathcal M$. 
By~\cite{bg-1}, we can choose 
 $\{\Omega_i\}$ as $\{\pa_i\Omega\}$, or $\{\nabla_i\Omega\}$. The latter are not holomorphic  on $\underline{\mathcal F}^{n-1}$ but are  holomorphic  on  $\underline{\mathcal H}^{n-1,1}$. We claim that
\begin{equation}\label{iol}
\bar\pa\nabla_i\Omega=g_{i\bar j}d\bar t_j\Omega.
\end{equation}
To see this, we first observe that
\[
\nabla_i\Omega=\pa_i\Omega-K_i\Omega,
\]
where $K_i=\pa_i\log(\Omega,\bar\Omega)$ are local functions. The choices of local functions $K_i$ are made so that $\nabla_i\Omega$ is orthogonal to  $\Omega$. Thus~\eqref{iol} follows from Theorem~\ref{thm61}. By the curvature formula~\eqref{1+2}, we have
\begin{equation}\label{1+2+3}
(R')_{i\bar j\gamma\bar\delta}
=(\ii)^{n-2}\left((\nabla_\gamma\nabla_i\Omega,
\overline{ \nabla_\delta\nabla_j\Omega})
-g_{i\bar\delta}g_{\gamma\bar j}(\Omega,\bar\Omega)\right).
\end{equation}

Using the above result, we get the following
\begin{lemma} \label{lem63}
Under the Kodaira-Spencer identification, we have
\[
R=I\omega_{WP}+(\ii)^{n}R',
\]
where $I$ is the identity map in $Hom(T\mathcal M,T\mathcal M)$.
\end{lemma}

{\bf Proof.} 
By Theorem~\ref{thm61}, we know that
\[
g_{i\bar j}=-\pa_i\bar\pa_j\log(\Omega,\bar\Omega).
\]
By a straightforward computation, we have
\[
g_{i\bar j}=-\frac{(\nabla_i\Omega,\bar{\nabla_j\Omega})}{(\Omega,\bar\Omega)}.
\]
The above formula gives the relation between the Weil-Petersson metric and the Hodge metric on  $\underline{\mathcal H}^{n-1,1}$. 
If we choose frames such that at a point, $(\sqrt{-1})^n(\Omega,\bar\Omega)=1$ and $g_{i\bar j}=\delta_{ij}$, then the metric of $\underline {\mathcal H}^{n-1,1}$ is also the identity  matrix at that point.
The lemma thus follows from~\eqref{1+2+3} and Theorem~\ref{thm62-1}.

\qed

The main result of this section is the following

\begin{theorem}\label{thm63}
Let $\gamma$ be a rational number, and let $f$ be an invariant polynomial on $Hom(T\mathcal M,T\mathcal M)$ with rational coefficients. Let $R_{WP}$ be the curvature tensor of the Weil-Petersson metric. Then we have
\begin{equation}\label{pklm}
\int_{\mathcal M} f(R_{WP}+\gamma I\omega_{WP})\in\mathbb Q,
\end{equation}
where $I$ is the identity map on $Hom(T\mathcal M,T\mathcal M)$. Moreover, we have
\begin{equation}\label{uio}
|c_{\alpha_1}(\omega_{WP})\wedge\cdots c_{\alpha_r}(\omega_{WP})\wedge\omega_{WP}^{\alpha_0}|
\leq 2^{m-\alpha_0}\,\omega_H^m
\end{equation}
for $\sum_{j=0}^r\alpha_j=m$, where $\omega_H=\omega_{H^n}$ is the Hodge metric.
\end{theorem}

{\bf Proof.} The equation~\eqref{pklm} follows from Theorem~\ref{thm41} and Lemma~\ref{lem63}.
The equation~\eqref{uio} essentially follows from Theorem~\ref{thmm52}. In what follows, we directly prove~\eqref{uio} in order to get the desired constant.

 First we choose a
 coordinate system at $x_0\in\mathcal  M$ so that  
\begin{equation}\label{dg}
g_{i\bar{j}}(x_0)=\delta_{ij}.
\end{equation}

Let $R=(R_i^j)$, and let 
\[ R_i^j = \sum_{k,l}R_{ik\bar{l}}^jdt_k\wedge d\bar{t}_l, \]
where 
$R^j_{ik\bar{l}}=g^{j\bar{p}}R_{i\bar{p}k\bar{l}}$.
Then the $\alpha$-th Chern class is given by
\begin{equation}\label{chern}
 c_\alpha(\omega_{WP})=\left(\bb\right)^\alpha\frac{(-1)^\alpha}{\alpha!}\sum_{\tau\in S_\alpha} sgn(\tau) R_{i_1}^{i_{\tau(1)}}\wedge \cdots \wedge R_{i_\alpha}^{i_{\tau(\alpha)}}, 
 \end{equation}
where $S_\alpha$ is the symmetric group on the set $\{1,2,\ldots,\alpha\}$.

We define

\[ h'_{i\bar j} = \gd_{ij}+\sum_l\langle \nabla_i \nabla_l\Omega,\bar{\nabla_j \nabla_l\Omega}\rangle. \]
Then $(h_{i\bar j}')$ defines a \ka metric $\omega'$.\footnote{The metric is equivalent to the partial Hodge metric in~\cite{ls-1}*{\S 4}. But we don't need this fact here.}

By Proposition~\ref{prop32}, we have
\[
\omega'\leq\omega_H.
\]
Thus in order to prove the theorem, we  only need to prove that
\[
 |c_{\alpha_1}(\omega_{WP})\wedge\cdots c_{\alpha_r}(\omega_{WP})\wedge\omega_{WP}^{\alpha_0}|
\leq 2^{m-\alpha_0}\,(\omega')^m.
\]

Let $A_{ij}=\sum_k(\nabla_i \nabla _k\Omega,\bar{\nabla _j\nabla _k\Omega})$. 
 Since the matrix $(A_{ij})$ is Hermitian, after a suitable
unitary change of basis, we can assume
\[ A_{ij}(x_0) = \left\{ \begin{array}{ll}
\lambda_i & \textrm{if $i=j$} \\
0 & \textrm{if $i\ne j$}.
\end{array}\right. \]		
Since $(A_{ij}(x_0))$ is semi-positive definite, $\lambda_i \ge 0$, and we have
\begin{equation}\label{dh}
h'_{i\bar{j}}(x_0) = \delta_{ij}(1+\lambda_i).
\end{equation}

For fixed $i,j,k,l$, by the Cauchy inequality, we have 
\begin{eqnarray*}
\big|R^j_{ik\bar{l}}\big| &\le& |\gd_{ij}\gd_{kl}+\gd_{il}\gd_{kj}- (\nabla _i\nabla _k\Omega,\bar{\nabla _j\nabla _l\Omega})| \\
&\le& 2+ \sqrt{\langle\nabla _i\nabla _k\Omega,\bar{\nabla _i\nabla _k\Omega}\rangle\langle\nabla _j\nabla _l\Omega,\bar{\nabla _j\nabla _l\Omega}\rangle}  \\
&\le& 2+ \sqrt{\lambda_k\lambda_l}\leq 2  \sqrt{(1+\lambda_k)(1+\lambda_l)}.
\end{eqnarray*}
So we get
\[ \Big|R_{j_1k_1\bar{l_1}}^{i_1}\ldots R_{j_\alpha k_\alpha \bar{l_\alpha}}^{i_\alpha} \Big| \le
2^\alpha\prod_{j=1}^\alpha \Big(\sqrt{(1+\lambda_{k_j})(1+\lambda_{l_j}}) \; \Big). \]
We note that all $k_j$'s (and also $l_j$'s)  have to be different. Thus we have
\[
\prod_j\sqrt{(1+\lambda_{k_j})(1+\lambda_{l_j})}\leq\det (h_{\alpha\bar\beta}).
\]
Using the same method as in the proof of Theorem~\ref{thmm52}, we get~\eqref{uio}.

\qed

\begin{cor}[Chern number inequalities]
Let $\alpha_0,\cdots,\alpha_r\geq 0$ be nonnegative numbers such that $\sum_{j=0}^r r_j=m$. Then
\[
c_{\alpha_1}(\omega_{WP})\wedge\cdots c_{\alpha_r}(\omega_{WP})\wedge \omega_{WP}^{\alpha_0}
\]
is absolutely integrable. Moreover, we have the following Chern number inequalities
\[
\int_M c_{\alpha_1}(\omega_{WP})\wedge\cdots c_{\alpha_r}(\omega_{WP})\wedge \omega_{WP}^{\alpha_0}
\leq 2^{m-\alpha_0}\int_M\omega^m_H.
\]
\end{cor}

{\bf Proof.} Since the Hodge volume is finite by~\cite{ls-1}*{Theorem 5.2}, the corollary follows.

\qed

  Theorem~\ref{thm63} is the counterpart of Theorem~\ref{thm41} and
  Theorem~\ref{thmm52} in the Weil-Petersson geometry. In what follows we shall show that
for Calabi-Yau moduli, we can do better than Hodge theory  provides. Since for application,  moduli space of a Calabi-Yau threefold is the most interesting one, for the rest of the paper, we will assume that  $\mathcal M$ is the Calabi-Yau moduli  of  a Calabi-Yau threefold.
 
   Recall that for the moduli space of a Calabi-Yau threefold, there is a simple relation between the Hodge metric and the
  Weil-Petersson metric~\cite{Lu5}:
  \begin{equation}\label{good}
  \omega_H=(m+3)\omega_{WP}+{\rm Ric}(\omega_{WP}).
  \end{equation}
 By the computation of~\cite{lu12}, since $\omega_H$ should be of Poincar\'e like at  infinity,
  it would be natural to conjecture the following
  
  \begin{conj} \label{conj1}The curvature of the Hodge metric $\omega_H$ is bounded.
  \end{conj}
  
  In this direction, 
we can prove the following
  \begin{theorem}\label{thmuj}
  The curvature of the Hodge metric is $L^1$ with respect to the Hodge metric.
  \end{theorem}

  \begin{remark} In ~\cites{Lu3,Lu5}, we have proved that the Ricci curvature of the Hodge metric is bounded from above by a negative number. Thus Theorem~\ref{thmuj} implies that the Hodge volume is finite, a result proved in~\cite{ls-1}. On the other hand, a slight modification of our proof also gives the Gauss-Bonnet-Chern theorm for the first Chern class of the Hodge metric.
  \[
  \int_{\mathcal M} c_1(\omega_H)\wedge\omega^{m-1}_H\in\mathbb Q.
  \]
  \end{remark}

  Before proving the theorem, we first give the following technical
  
  \begin{lemma}\label{lem62} Let $U$ be a small neighborhood of $\C^m$ at the origin. Assume that $(s_1,\cdots,s_a,w_1,\cdots,w_b)$ are the local coordinates ($a+b=m$) similar to those   in Lemma~\ref{lem42}. Let $\omega_P$ be the Poincar\'e metric
  \[
\omega_P=\bb\left(\sum_{i=1}^a\frac{ds_i\wedge d\bar s_i}{|s_i|^2(\log\frac{1}{|s_i|})^2}
+\sum_{i=a+1}^{a+b}dw_i\wedge d\bar w_i\right).
\] 
  Let 
  \[
  \phi_\eps(s_1,\cdots,s_{a})=\prod_{j=1}^{a} (1-\phi_\eps (s_j)),
  \]
  where $\phi_\eps$ is defined in ~\eqref{phi-1}.
  Then we have
  \[
  \int_{U\cap \,{\rm supp}\, \phi_\eps} \left(\sum_{j=1}^a\log|s_j|\right)\omega_P^{m}\leq C,
  \]
  where $C$ is independent to $\eps>0$.
  \end{lemma}
  
  {\bf Proof.} We observe that the support of $\pa\bar\pa\phi_\eps$ is the union
$$\bigcup_{j=1}^a \{e^{-\frac 1\eps}\leq |s_j|\leq e^{-\frac{1}{2\eps}}\}.$$ 
We assume that $U\subset \{|s_j|\leq \frac 12, |w_j|\leq \frac 12\}$.
Then the  expression is less than a constant times
\[
\sum_{j=1}^a\left(\int_{e^{-\frac 1\eps}}^{e^{-\frac 1{2\eps}}}\frac{1}{|s_j|\log\frac{1}{|s_j|}} d|s_j|\prod_{i\neq j}\int_0^{\frac 12}\frac{1}{|s_i|(\log\frac{1}{|s_i|})^2} d|s_i|\right)\leq C.
\]
\qed

{\bf Proof of Theorem~\ref{thmuj}.} 
By~\cite{Lu5}, we know that the norm of the curvature is bounded by the Ricci curvature of the Hodge metric. Thus we just need to prove that
\[
\int_{\mathcal M} c_1(\omega_H)\wedge\omega^{m-1}_H\leq C.
\]

The proof depends on the analysis of the curvature and the Hodge metric itself at infinity. Thus as before,  we assume that $M$ can be compactified to be a compact manifold $\bar M$ by a divisor $Y$ of normal crossings.
  
  By~\eqref{good} and Theorem~\ref{thm62-1}, we know that $|c_1(\omega_{WP})|\leq\omega_H$. Thus we have
\[
\int_{\mathcal M} c_1(\omega_{WP})\wedge\omega^{m-1}_H\leq C
\]
by Corollary~\ref{cor53}.
Combining the above two inequalities, we need to prove that 
\[
\bb\int_{\mathcal M} \pa\bar\pa\log\frac{\omega_H^m}{\omega_{WP}^m}\wedge\omega^{m-1}_H\leq C.
\]
Let $\rho_\eps$ be the cut-off function defined in Lemma~\ref{lem42}. Using integration by parts, we concluded that the following inequality implies Theorem~\ref{thmuj}:
\[
\bb\int_{\mathcal M} \log\frac{\omega_H^m}{\omega_{WP}^m}\wedge\pa\bar\pa\rho_\eps\wedge\omega^{m-1}_H\leq C,
  \]
 where $C$ is independent to $\eps$.
  
  Let $V_\eps$ be the support of the form $\pa\bar\pa\rho_\eps$. By the properties of $\rho_\eps$,  up to a constant, we have
  \[
 \bb \int_{\mathcal M} \log\frac{\omega_H^m}{\omega_{WP}^m}\wedge\pa\bar\pa\rho_\eps\wedge\omega^{m-1}_H\leq
  \int_{V_\eps} \log\frac{\omega_H^m}{\omega_{WP}^m}\wedge\omega_P\wedge\omega_H^{m-1}.
    \]
  By Corollary~\ref{cor53} again, we know that up to a constant, we have
  \[
  \omega_H\leq \omega_P.
  \]
  Thus in order to  prove the theorem, we only need to prove that
  \[
\int_{V_\eps}\log\frac{\omega^m_P}{\omega^m_{WP}}\wedge\omega_P^m\leq C.
  \]
The problem being local, we may consider a neighborhood $U$ in Lemma~\ref{lem62}. By ~\cite{ls-2}*{Lemma 6.4}, there is a local function $f$ such that
  \[
\omega_{WP}^m=\left(\prod_{j=1}^a |s_j|^{p_j}\right)\,f
ds_1\wedge\cdots d\bar s_a\wedge dw_1\wedge\cdots\wedge d\bar w_b,
  \]
  and $\log f$ is integrable with respect to the metric $\omega_P$. By the above equation, we have
  \[
  \frac{\omega_P^m}{\omega_{WP}^m}=\frac{\prod_{j=1}^a (\log\frac{1}{|s_j|})^2}{
  \prod_{j=1}^a |s_j|^{p_j+2}\,f}.
    \]
    Since $\log\log \frac{1}{|s_j|}$  is integrable with respect to $\omega_P$, the theorem follows from Lemma~\ref{lem62} and the fact that $\log f$ is also  integrable with respect to $\omega_P$.
    
    \qed
    
\begin{remark} Theorem~\ref{thmuj} is a generalization of the fact that the Hodge metrics are Poincar\'e bounded. 
At the moment, it seems less well motivated. However, we suspect that for the Calabi-Yau moduli, we can prove one more layer of estimate provided by the variation of Hodge structure. 
Such an estimate will reflect the inner structure of a Calabi-Yau manifold (cf.~\cite{lu12} and Conjecture~\ref{conj1}).  Thus this theorem may point
the way for future research.
\end{remark}

We end our paper with the following physics applications. First, by
Theorem~\ref{thmuj}, the volume of the Weil-Petersson metric is
finite (cf.~\cite{ls-1}). The result has the following string theoretical
implications (cf.~\cite{dl-1}).

In ~\cite{AD} and ~\cite{DSZ-3}, the index of all supersymmetric vacua
was given. In~\cite{AD}*{eq. (1.5)}, the index is given by
  \[
  I_{vac}(L\leq L_{max})=const. \int_{\mathcal M\times\mathcal H}\det(-R_{WP}-\omega_{WP}),
  \]
  where $\mathcal M$ is the Calabi-Yau moduli and $\mathcal H$ is the moduli space of elliptic curves. In~\cite{DSZ-3}*{Theorem 1.8}, The following strengthened  result of the above was given.
  
  \begin{theorem} Let $K$ be a compact subset of $\mathcal M$ with piecewise smooth boundary. Then
  \[
  {\mathcal I}nd_{\chi_K}(L)=const. (L^{2m})\left[\int_K c_m(T^*(M)\otimes\underline{\mathcal F}^3)+O(L^{-1/2})\right].
  \]
  \end{theorem}
  
  By Theorem~\ref{thm63}, we have
  \begin{theorem} The indices $I_{vac}$ and ${\mathcal I}nd_{\chi_K}$ are all finite. Moreover, 
  ${\mathcal I}nd_{\chi_K}$ is bounded from above uniformly with respect to $K$.
  They are all bounded, up to an absolute constant, by the Hodge volume of the Calabi-Yau moduli.
  \end{theorem}
  
  \qed

\end{document}